\documentclass[final,12pt]{article}

\usepackage[final]{graphics}

\usepackage{amsmath,amsfonts,amsthm,upref}

\usepackage[fracintave,numstogether]{mathdef}

\usepackage{mylabels}

\usepackage{mymargs}

\usepackage{verbatim}

\begin{document}
\pagestyle{myheadings}
\sloppy


\author{Harri Ojanen%
\thanks{The author was supported in part by the Academy of Finland and
    the Finnish Society of Sciences and Letters.}}
    
\title{\bf Weighted Norm Estimates and\\ Representation Formulas for\\ Rough
	Singular Integrals}

\date{June 30, 1998}

\maketitle

{\footnotesize
\begin{quote}
{\bf Abstract}:
Weighted norm estimates and representation formulas are proved for
non-homogeneous singular integrals with no regularity condition on the
kernel and only an $L\log L$ integrability condition. The representation
formulas involve averages over a star-shaped set naturally associated
with the kernel. The proof of the norm estimates is based on the
representation formulas, some new variations of the Hardy-Littlewood
maximal function, and weighted Littlewood-Paley theory.

\medbreak



\end{quote}
}


\markright{Rough Singular Integrals}
\thispagestyle{empty}


\section*{Introduction}

R.~Fefferman~\cite{Feff79} introduced non-homogeneous singular integral
operators of the form
\[
    T f(x) = \pv \int_{\Rn}\frac{\Omega(y)}{\abs{y}^n}
    h(\abs{y}) f(x-y)\, dy.
\]
We prove conditions on weight functions so that $T$~is a bounded mapping
on the weighted space $\Lpw = L^p(\Rn,w(x)dx)$, $1<p<\infty$, $n\geq2$. 
In addition to the standard requirements on the kernel
($\Omega$~is positively homogeneous of degree zero and $\int_\Sn
\Omega(\theta)\,d\theta=\nolinebreak0$),
we assume only that $\Omega$~is in $\LlogLS$, that is,
\[
    \snorm{\Omega}{\LlogL} = \int_\Sn \abs{\Omega(\theta)}
                (1+\logp\abs{\Omega(\theta)}) \, d\theta < \infty.
\]
In particular, there is no regularity condition on~$\Omega$.
The radial function~$h$ satisfies
\[
    \int_R^{2R} \abs{h(r)}^\sigma\,dr \leq C_hR
\]
for all $R>0$ and for an appropriate value of $\sigma>1$.

J.~Duoandikoetxea and J.L.~Rubio de Francia~\cite{Duoa86} proved
weighted norm inequalities for~$T$ when~$\Omega$ is essentially bounded
on~$\Sn$. Their results were generalized by D.K.~Watson~\cite{Wats90} to
the case $\Omega\in L^r(\Sn)$, for some $1<r\leq\infty$. We extend the
results of \cite{Duoa86} and~\cite{Wats90} by considering a more general
class of kernels, namely $\Omega$ in $\LlogLS$. Our results improve
those proved by D.K.~Watson and R.L.~Wheeden~\cite{Wats97}, who studied
homogeneous operators, i.e., the case $h\equiv1$, with
$\Omega\in\LlogLS$. We are also able to prove weighted inequalities and
representation formulas for truncated singular integrals, which were not
obtained in~\cite{Wats97}.

The weighted norm estimates are based on a representation formula for
the singular integral using averages over a star-shaped set that is naturally
associated with the operator. The set consists of those points where
$\abs{\Omega(x)}/\abs{x}^n$ is greater than~$1$ and is denoted by~$S_\Omega$.
The homogeneity of~$\Omega$ implies
$S_\Omega$~is star-shaped about the origin, i.e., if $x\in S_\Omega$ then $tx\in
S_\Omega$ for all $0<t<1$. We call $S_\Omega$~the star-shaped set
associated with~$\Omega$.

The approach used in~\cite{Duoa86} and~\cite{Wats90}, when $\Omega\in
L^r(\Sn)$, $1<r\leq\infty$, is to study weighted estimates for $T$ by
using the Muckenhoupt $A_p$~weights. This class consists of those
positive locally integrable functions~$w$ for which
\begin{equation}
    \sup_Q
    \biggl( \intave{Q}{w(x)\,dx} \biggr)^{1/p}
    \biggl( \intave{Q}{w(x)^{-p'/p}\,dx} \biggr)^{1/p'}   \tag{$A_p$}
	< \infty,
\end{equation}
the supremum is taken over all cubes $Q\subset \Rn$
($p'$~is the dual to $p$ defined by $1/p+1/p'=1$).

Since we assume that the homogeneous part~$\Omega$ of the kernel lies in
$\LlogLS$ and not necessarily in any $L^r(\Sn)$ for $r>1$, the structure
of the set~$S_\Omega$ yields restrictions on the weight: we require~$w$
to satisfy a condition similar to an $A_p$~condition, but with
\textit{rectangles} related to the set~$S_\Omega$ instead of
cubes. In general this is a more restrictive condition than the
$A_p$~condition, since the eccentricities of these rectangles may be
unbounded. See Theorems~\ref{t:Thbded}, \ref{t:Thinfbded}, and
Corollary~\ref{t:Thbdedunif} for details. (The approach of using
star-shaped sets is interesting also when $\Omega\in L^r(\Sn)$: the
results of~\cite{Wats90}, in the special case $h\equiv1$, are derived
in~\cite{Wats97} by using this method.)

The results are based on a representation formula for truncated singular
integrals
\[
    T_{\epsilon}f(x) = \int_{\abs{y}>\epsilon}\frac{\Omega(y)}{\abs{y}^n}
    h(\abs{y}) f(x-y)\, dy.
\]
In Theorem \ref{t:srep} we show that these operators
can be written in terms of averages over dilates of
the set~$S_\Omega$: in fact,
$
    T_{\epsilon}f(x) = n \int_0^\infty A_{\epsilon,t} f(x) \,dt/t,
$
where $A_{\epsilon,t}$ is the ``average''
\[
    A_{\epsilon,t}f(x) = \frac1{t^n}
    \int_{tS_\Omega\setminus B(0,\epsilon)}
    f(x-y) h(\abs{y}) \sgn\Omega(y)\,dy.
\]
See Theorem \ref{t:srep} for the exact statement.  We
also show similar formulas for the principal value operator~$T$ and for
some classes of non-convolution type operators.
As an application we discuss the Calder\'on commutators.

The weighted norm estimates for $T$ and $T_\epsilon$
require the study of an associated maximal operator defined by
\[
    M_h f(x) = \sup_{r>0} \frac1{r^n}
        \int_{B(x,r)}{\mabs{h(\abs{x-y}) f(y)}\,dy}.
\]
We prove weighted norm estimates and vector-valued inequalities for
slightly more general maximal operators as well as weighted norm
estimates for a related ``starlike operator'' (where the integration is
over dilates of a star-shaped set instead of balls). See Theorems
\ref{t:maxApr}, \ref{t:maxvect}, and \ref{t:MSh}, respectively. We
also discuss results for corresponding fractional maximal operators.

The content of each section is as follows: Section~\ref{s:res} contains
the statements of the theorems on singular integrals. The maximal
operators are studied in Section~\ref{s:max}. Section~\ref{s:srep}
contains the proof of the representation formula for truncated
operators. Some preliminary results are given in Section~\ref{s:prelim}
before the proof of the weighted norm inequalities for singular integrals
in Section~\ref{s:hbded}. Finally in Appendix~\ref{s:hrep} we state
and prove representation formulas for non-convolution type and principal
value operators.

The proof of Theorem~\ref{t:Thbded} is based on similar
arguments in~\cite{Wats97} for homogeneous singular integrals.  I thank
professor Richard Wheeden for his suggestions and encouragement during
this project and professor David Watson for useful comments.

As usual the letter~$C$ denotes a constant whose value may change from
one line to the next.


\section{Statement of main results}\mylabel{s:res}

\begin{definition}\mylabel{d:S}
Let $\Omega\in L^1(\Sn)$, $n\geq2$, be positively homogeneous of
degree~$0$. The
\textit{star-shaped set $S_\Omega \subset \Rn$ associated with $\Omega$}
is 
$
    S_\Omega = \{x\in\Rn \colon \abs{x} \leq \rho_\Omega(x) \},
$
where $\rho_\Omega(x) = \abs{\Omega(x)}^{1/n}$.
\end{definition}

\begin{definition}\mylabel{d:H}
We let $\setH(\sigma)$, $1\leq\sigma<\infty$,
denote the class of all measurable complex valued
functions defined on $\Rp$ that satisfy the following condition: 
there exists $C_h>0$ such that
$\int_R^{2R} \abs{h(r)}^\sigma\,dr \leq C_hR$ for all $R>0$.
\end{definition}

Note that by H\"older's inequality $\setH(\sigma_1)$ is contained in
$\setH(\sigma_2)$ when $\sigma_1>\sigma_2$.

\begin{theorem}\mylabel{t:srep}
Let $\Omega\in L^1(\Sn)$ be positively homogeneous of degree~$0$, 
$h\in\setH(p')$,
and $f\in L^p(\Rn)$, $1<p<\infty$. For $\epsilon>0$ define the operator
\[
    T_{\epsilon}f(x) = \int_{\abs{y}>\epsilon}\frac{\Omega(y)}{\abs{y}^n}
    h(\abs{y}) f(x-y)\, dy
\]
and let
\begin{subequations}\mylabel{e:srep}
\begin{equation}\mylabel{e:srepX}
    A_{\epsilon,t} f(x) = \frac1{t^n} \int_{tS\setminus 
    B(0,\epsilon)} f(x-y) h(\abs{y})\sgn\Omega(y)\,dy,
\end{equation}
where $S=S_\Omega$ is the star-shaped set associated with $\Omega$. Then
for almost all $x\in\Rn$ the representation formula
\begin{equation}\mylabel{e:srepY}
    T_\epsilon f(x) = n \int_0^\infty A_{\epsilon,t} f(x) \frac{dt}t
\end{equation}
\end{subequations}
holds 
and the integrals in \eqref{e:srepX} and \eqref{e:srepY} converge 
absolutely.

Moreover, if $f\in\Li(\Rn)$ has compact support and $h\in\setH(1)$ the
representation formula \eqref{e:srepY} holds for all $x\in\Rn$
and the integrals in \eqref{e:srepX} and \eqref{e:srepY} converge 
absolutely for all $x\in\Rn$.
\end{theorem}

See also Appendix~\ref{s:hrep} for a discussion of representation
formulas for principal value and non-convolution type operators.  Note
that in~\cite{Wats97} the authors use a representation formula for
principal value operators, i.e., a result more like
theorem~\ref{t:hrep} than~\ref{t:srep}. The nature of the
operators studied here makes it more convenient to use a representation
formula for truncated operators.

\begin{remark}\mylabel{r:limTheps}
By writing the truncated operators in polar coordinates it is easy to
show that $\lim_{\epsilon\tozerop} T_\epsilon f(x)$ exists for all
$x\in\Rn$ on test functions, say when $f\in C_0^1(\Rn)$ and $h\in\setH(1)$.
Moreover, the convergence is uniform in $\Rn$.
The key observation is $\lim_{\epsilon\tozerop} \int_0^\epsilon
\abs{h(r)}\,dr=0$, see Lemma \ref{l:Heq}(b). In fact, for
$\epsilon>\eta>0$, we have
\[
    \abs{T_\eta f(x) - T_\epsilon f(x)}
    \leq
    \int_\Sn \abs{\Omega(\theta)}
        \int_\eta^\epsilon \mabs{h(r)(f(x-r\theta)-f(x))}
            \frac{dr}r \,d\theta,
\]
which is bounded by $\snorm{\Omega}{1} \snorm{\nabla f}{\infty}
\int_0^\epsilon \abs{h(r)}\,dr$. Letting $\epsilon\tozerop$
and using Lemma \ref{l:Heq}(b) shows the principal
value operator $T$ is well defined as the pointwise limit
at least on $C_0^1(\Rn)$.
\end{remark}

In the following the set $S = S_\Omega$ is always the star-shaped set
associated with $\Omega$ and $\rho=\rho_\Omega$ (see definition
\ref{d:S}). We decompose $S$ as a disjoint union $S=\bigcup_{m=0}^\infty S_m$,
where
\begin{align}\mylabel{e:Smdef}
    S_m &= \{ x\in S\colon 2^{m-1}<\rho(x)\leq2^m \}
        \text{, for $m\geq1$, and } \\
    S_0 &= \{ x\in S\colon \rho(x)\leq1 \}, \notag
\end{align}
and we also use 
the corresponding projections on the unit sphere
\begin{align}\mylabel{e:Thmdef}
    \Theta_m &= \{ \theta\in\Sn\colon 2^{m-1}<\rho(\theta)\leq2^m \}
        \text{, for $m\geq1$, and } \\
    \Theta_0 &= \{ \theta\in\Sn\colon \rho(\theta)\leq1 \}. \notag
\end{align}

\begin{definition}\mylabel{d:starcov}
Let $S$ be as above. A \textit{stratified starlike cover} of $S$ is a
collection of rectangles $\{R_{m,k}\}$,  $m\geq0$, $1\leq k<k_m$ with
$0\leq k_m\leq\infty$, that satisfies the
following conditions:
\begin{enumerate}
\item Each $R_{m,k}$ is centered at the origin and  the length of the
    longest side of $R_{m,k}$ is comparable to $2^m$.
\item For all $m\geq0$, 
    $S_m\subset 
    \bigcup_k R_{m,k}$ up to a set of measure zero and
    $\sum_k \abs{R_{m,k}} \leq c_n \abs{S_m}$, where $c_n$ depends
    only on the dimension $n$.
\end{enumerate}
\end{definition}

\begin{remark}\mylabel{r:coverconstr}
It is easy to show that such a cover always exists, see
\cite{Wats97} or \cite[pp.~248--249]{Chan93}.
The idea is to cover each $\Theta_m\subset\Sn$ by ``disks''
$D_{m,k}\subset\Sn$, whose combined (surface) measure is proportional to
the measure of $\Theta_m$. The disks are used to define cones $C_{m,k}$ with
vertex at the origin, intersection with $\Sn$ equal to $D_{m,k}$, and
height~$2^m$. Clearly $\bigcup_k C_{m,k}$ covers $S_m$ up to a set of measure
zero. Then $R_{m,k}$ is defined to be any one of the smallest rectangles
centered at the origin that contain~$C_{m,k}$.

A very simple two-dimensional example is shown in figure~\ref{f:sets2},
where, for some $\alpha\in(0,1)$, $\Omega$ satisfies $\abs{\Omega(\theta)} =
\abs{\sin\theta}^{-\alpha}$, $0\leq\theta<2\pi$.
The set $S_\Omega$ has two unbounded ``arms'' along the $x_1$-axis. The
cover shown is of the form $\{R_{m,k}\}$,
$m\geq0$, $k=1,2$. The rectangles become wider and more eccentric and their
major axis turns towards the $x_1$-axis as $m$ increases.
\begin{figure}
\begin{center}
\scalebox{0.75}{\includegraphics{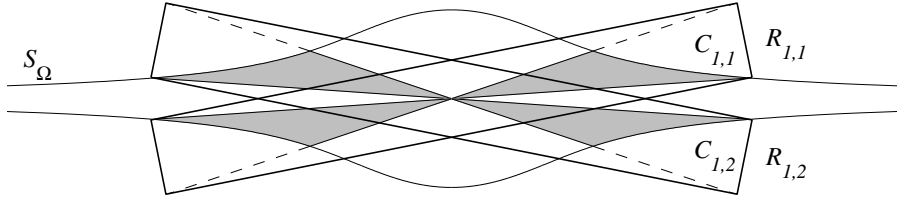}}
\end{center}
\caption{Part of a stratified starlike cover of $S_\Omega$: 
    the shaded region, $S_1$, is
    covered by rectangles $R_{1,1}$ and $R_{1,2}$. 
    See remark \ref{r:coverconstr}.
    }
\mylabel{f:sets2}
\end{figure}
\end{remark}

For any set $E\subset\Rn$ the collection $\B(E)$ consists of all 
translates and (isotropic) dilates of $E$.

\begin{allowdisplaybreaks}
\begin{theorem}\mylabel{t:Thbded}
Suppose $\Omega\in\LlogLS$ is positively homogeneous of degree~$0$ and
$\int_\Sn\Omega(\theta)\,d\theta=0$. Let
$\{R_{m,k}\colon m\geq0, 1\leq k<k_m\}$ be a stratified starlike cover
of the set $S=S_\Omega$.
Assume $1<p<\infty$ and that $w$ is a non-negative measurable
function that satisfies the following condition:
    there exists $r>1$ such that for all $m\geq0$, $1\leq k<k_m$
    there is a constant $K_{m,k}>0$ so that
    \begin{subequations}\mylabel{e:ThcondC}
    \begin{align}
        \mylabel{e:ThcondCa}
        \biggl(\intave{R}{w}\biggr)^{1/p}
        \biggl(\intave{R}{w^{-rp'/p}}\biggr)^{1/rp'}
        &\leq \frac{K_{m,k}}{\abs{R_{m,k}}},
        && \text{if $1<p\leq2$,}
        \\
        \mylabel{e:ThcondCb}    
        \biggl(\intave{R}{w^r}\biggr)^{1/rp}
        \biggl(\intave{R}{w^{-p'/p}}\biggr)^{1/p'}
        &\leq \frac{K_{m,k}}{\abs{R_{m,k}}},
        && \text{if $2\leq p<\infty$,}
    \end{align}
    \end{subequations}
    holds for all $R\in\B(R_{m,k})$  and the constants $K_{m,k}$ satisfy 
    \begin{equation}\mylabel{e:Thsumcond}
        \sum_{m=0}^\infty \sum_{k=1}^{k_m-1} (m+1) K_{m,k} < \infty.
    \end{equation}
Then there exists $\sigma>1$ such that if
$h\in\setH(\sigma)$ the operators $T_\epsilon$, $\epsilon>0$, can be
extended to a uniformly bounded family of operators on $\Lpw$ and the
operator~$T$
can be extended to a bounded operator on $\Lpw$.
\end{theorem}
\end{allowdisplaybreaks}

Theorem 1.5 of \cite{Wats97} is the corresponding result when
$h\equiv1$. It is included as a special case of the above theorem. More
generally, if $h\in\Li(\Rp)$ then clearly $h\in\setH(\sigma)$ for all
$\sigma>1$. Hence, in this case, whether $T$ is bounded does not depend
on the properties of $h$.

\begin{remark}\mylabel{r:wAp}
Condition \ref{e:ThcondC} implies $w\in A_p$:
Since \eqref{e:ThcondC} holds for all translates and
dilates of any of the fixed rectangles $R_{m,k}$,
it holds for all cubes in place of the
rectangles. Then H\"older's inequality shows $w\in A_p$.
\end{remark}

\begin{corollary}\mylabel{t:Thbdedunif}
Suppose $\Omega$ and the sets $\{R_{m,k}\}$ are as in
Theorem~\textup{\ref{t:Thbded}},
and that instead of \eqref{e:ThcondC}
    there exists a constant $K>0$ such that for all $m\geq0$,
    $1\leq k<k_m$ and $R\in\B(R_{m,k})$
    \begin{equation}\mylabel{e:ThcondA}
        \biggl(\intave{R}{w}\biggr)^{1/p}
        \biggl(\intave{R}{w^{-p'/p}}\biggr)^{1/p'} \leq K.
    \end{equation} 
Then there exists $\sigma>1$ such that if
$h\in\setH(\sigma)$ then $T$ can be extended to a bounded operator on
$\Lpw$ and the truncated operators $T_\epsilon$, $\epsilon>0$, can be extended to
a uniformly bounded family of operators on $\Lpw$. \end{corollary}

\begin{proof}
Since \eqref{e:ThcondA} is uniform in $m$ and $k$
the reverse H\"older's inequality applied to $w$ and $w^{-p'/p}$ shows
that $w$ satisfies condition \eqref{e:ThcondC} with $K_{m,k}=K\abs{R_{m,k}}$ and
for some $r>1$, see \cite{Wats97}.
Then \eqref{e:Thsumcond} is a consequence of the hypothesis $\Omega\in\LlogLS$
and of $\sum_k \abs{R_{m,k}} \leq c \abs{S_m}$,
see \eqref{e:mSmsum}.
\end{proof}

\begin{theorem}\mylabel{t:Thinfbded}
Suppose $\Omega$ and the sets $\{R_{m,k}\}$ are as in
Theorem~\textup{\ref{t:Thbded}},
$h\in\Li(\Rp)$, and for all $m\geq0$, $1\leq k<k_m$
there exists a constant $K_{m,k}>0$ such that for all $R\in\B(R_{m,k})$
\begin{equation}\mylabel{e:ThcondB}
    \biggl(\intave{R}{w}\biggr)^{1/2}
    \biggl(\intave{R}{w^{-1}}\biggr)^{1/2}
    \leq \frac{K_{m,k}}{\abs{R_{m,k}}},
\end{equation}
and the constants $K_{m,k}$ satisfy \eqref{e:Thsumcond}. 
Then 
$T$~can be extended
to a bounded operator on $L^2_w$,
and the
truncated operators $T_\epsilon$, $\epsilon>0$, can be extended to a
uniformly bounded family of operators on $L^2_w$.
\end{theorem}



\section{Maximal operators}\mylabel{s:max}

We define the maximal operator
\begin{equation}\mylabel{e:maxdef}
    M_H f(x) = \sup_{r>0} \frac1{r^n}
        \int_{\abs{y}<r}{H(x,y) \abs{f(x-y)}\,dy},
\end{equation}
where the non-negative function $H$ satisfies for some $\sigma>1$
and $C_H>0$
\begin{equation}\mylabel{e:hcube}
    \int_{\abs{y}<r}{ H(x,y)^\sigma\,dy} \leq C_H r^n
\end{equation}
for all $x\in\Rn$ and $r>0$.

There are two simple examples of functions that satisfy \eqref{e:hcube}.
One is given by radial functions in the class $\setH(\sigma)$, i.e.,
$H(x,y)=\abs{h_x(\abs{y})}$ and $h_x\in\setH(\sigma)$ uniformly in
$x\in\Rn$. Another example is given by homogeneous functions: let
$H(x,y)=\abs{\Phi_x(y)}$, where $\Phi_x$ is homogeneous of degree zero and
$\Phi_x\in L^\sigma(\Sn)$ uniformly in $x\in\Rn$. In both cases the
verification of \eqref{e:hcube} follows by writing the integral in
\eqref{e:hcube} in polar coordinates (see also Lemma~\ref{l:hstarint}
below).

In this section the letter $C$ denotes a constant that depends on $H$
only through the constant $C_H$ in equation \eqref{e:hcube} above.


\subsection{Weighted norm inequalities}

In order to study the starlike maximal function $M_{S,H}$ (see
Theorem~\ref{t:MSh} below) we need to know precisely how the operator norm of
$M_H$ on $\Lpw$ depends on the weight.

\begin{theorem}\mylabel{t:maxApr}
Let $1<p<\infty$ and assume that the weight $w$ satisfies for some
$r>1$
\begin{equation}\mylabel{e:Apr}
    \left(\intave{Q}{w}\right)^{1/p}
    \left(\intave{Q}{w^{-rp'/p}}\right)^{1/rp'} \leq K
\end{equation}
for all cubes $Q\subset\Rn$. If $H$ satisfies \eqref{e:hcube} for some
$\sigma>r'p'$ then $M_H$ is bounded on $\Lpw$ with
operator norm bounded by $CK$.
\end{theorem}

Theorem \ref{t:maxApr} follows easily from the next result, which
is a special case of Theorem 2.11 in~\cite{Pere94}. 
The Hardy-Littlewood maximal operator is defined by
$Mf(x) = \sup_{r>0} r^{-n} \int_{\abs{y}<r}{\abs{f(x-y)}\,dy}$.

\begin{theorem}[P\'erez]\mylabel{t:perez}
Let $1<p<\infty$ and suppose $w$ satisfies \eqref{e:Apr} for some
$r>1$. Then $\snorm{Mf}{p,w} \leq CK \snorm{f}{p,w}$, where $C$
is independent of $f$ and $K$.
\end{theorem}

The estimate $CK$ for the operator norm is not stated in~\cite{Pere94},
but it is easy to see that it follows from the proof given there.

\begin{proof}[Proof of Theorem \ref{t:maxApr}]
By H\"older's inequality we have
\[
    M_H f(x) \leq \sup_{r>0}
        \left( \frac1{r^n}\int_{\abs{y}<r}{H(x,y)^\sigma\,dy} \right)^{1/\sigma}
        \left( \frac1{r^n}\int_{\abs{y}<r}{\abs{f(x-y)}^{\sigma'}\,dy} \right)^{1/\sigma'}.
\]
Since $h$ satisfies~\eqref{e:hcube} we get
$M_Hf(x) \leq C_H^{1/\sigma} \bigl[M(\abs{f}^{\sigma'})(x)\bigr]^{1/\sigma'}$,
where $M$ is the Hardy-Littlewood maximal function. 
Therefore $\snorm{M_Hf}{p,w} \leq C
\smnorm[1/{\sigma'}]{M(\abs{f}^{\sigma'})}{p_1,w}$, where $p_1=p/\sigma'$.

Condition \eqref{e:Apr} implies a similar condition with the power~$p_1$
in place of~$p$: generalizing an idea from~\cite{Pere96} we claim that
there exists $\rho>1$ such that
\begin{equation}\mylabel{e:Aqrho}
    \left(\intave{Q}{w}\right)^{1/p_1}
    \left(\intave{Q}{w^{-\rho p_1'/p_1}}\right)^{1/\rho p_1'} \leq K^{\sigma'}
\end{equation}
for all cubes $Q\subset\Rn$. To prove this claim note that the equation
$rp'/p=\rho p_1'/p_1$ is equivalent with $\rho=r(p_1-1)/(p-1)$. It is easy to
see that $\rho>1$ is equivalent with $\sigma>r'p'$. 
Raising both sides of \eqref{e:Apr} to the power $\sigma'$ gives
\eqref{e:Aqrho}. Now Theorem~\ref{t:perez} applied on $L^{p_1}_w$ gives
$
    \snorm{M_Hf}{p,w} \leq C
    \bigl(K^{\sigma'} \smnorm{\abs{f}^{\sigma'}}{p_1,w}\bigr)^{1/\sigma'}
    = CK \snorm{f}{p,w}.
$
\end{proof}

\begin{corollary}\mylabel{c:maxAp}
Let $1<p<\infty$ and $w\in A_p$. Then there exists $\sigma>1$ such that
if $H$ satisfies \eqref{e:hcube} then $M_H$ is bounded on $\Lpw$.
\end{corollary}

\begin{proof}
The reverse H\"older's inequality implies $w$ satisfies \eqref{e:Apr} for
some $r>1$, see~\cite{Coif74b,Muck72}.
\end{proof}


\subsection{A starlike maximal operator}

If $S\subset\Rn$ is an arbitrary measurable star-shaped set centered at the
origin, a collection $\{R_j\}$ of rectangles is a \textit{starlike
cover} of $S$ if each $R_j$ is a rectangle centered at the origin,
$S\subset\bigcup_j R_j$ up to a set of measure zero, and
$\sum_j \abs{R_j} \leq c_n \abs{S}$,
where $c_n$ depends only on the dimension of $\Rn$. It
is shown in~\cite{Chan93} that such a cover always exists (see also
remark \ref{r:coverconstr}).

\begin{theorem}\mylabel{t:MSh}
Assume $H$ is a non-negative measurable function on $\Rn\times\Rn$,
$S\subset\Rn$ is star-shaped about the origin,
and $\{R_j\}$ is a
starlike cover of $S$. Define the starlike maximal operator
\begin{equation}\mylabel{e:MShdef}
	M_{S,H}f(x) = \sup_{t>0} \frac1{t^n} \int_{tS} H(x,y) \abs{f(x-y)} \, dy.
\end{equation}
Let $1<p<\infty$ and suppose the weight $w$ satisfies for some $r>1$
\begin{equation}\mylabel{e:AprRj}
    \left(\intave{R}{w}\right)^{1/p}
    \left(\intave{R}{w^{-rp'/p}}\right)^{1/rp'} \leq \frac{K_j}{\abs{R_j}},
\end{equation}
for all $R\in\B(R_j)$,
where the constants $K_j$ satisfy $\sum_j K_j < \infty$.
If $H$ satisfies 
\begin{equation}\mylabel{e:hrect}
    \int_{tR_j} H(x,y)^\sigma\,dy \leq C \abs{tR_j},
    \quad
    x\in\Rn,
\end{equation}
for all $t>0$ and all $j$, and for some $\sigma>r'p'$,
then $M_{S,H}$ is bounded on $\Lpw$ with operator norm bounded by
$C\sum_j K_j$, where $C$ is independent of $\{K_j\}$.
\end{theorem}


\begin{remark}
Theorem \ref{t:MSh} generalizes to a fractional version of the operator
$M_{S,H}$ defined by
\[
    M_{\mu,S,H}f(x) = \sup_{t>0} \frac1{t^{n-\mu}}
        \int_{tS} H(x,y) \abs{f(x-y)} \, dy, \quad 0\leq \mu<n.
\]
The result proved in~\cite{Pere94}, Theorem 2.11, is much more general
than Theorem~\ref{t:perez} stated above. Theorem 2.11 gives conditions
when the fractional maximal operator is bounded from $L^p_w$ to $L^q_v$,
$1<p\leq q<\infty$. The same proof that is given above for
Theorem~\ref{t:maxApr} combined with the full strength of Theorem 2.11
of~\cite{Pere94} can be easily used to get two weight norm estimates for
another fractional maximal operator,
\[
    M_{\mu,H} f(x) = \sup_{r>0} \frac1{r^{n-\mu}}
        \int_{\abs{y}<r} H(x,y) \abs{f(x-y)}\,dy, \quad 0\leq \mu<n,
\]
from $L^p_w$ to $L^q_v$, when the weights $v$ and
$w$ satisfy the same conditions as in~\cite{Pere94}, and
$\sigma>\max\{r'p', n/(n-\mu)\}$. Using these results Theorem
\ref{t:MSh} generalizes to $M_{\mu,S,H}$ 
when the weights $v$ and $w$ satisfy the conditions of Theorem 5(C)
of~\cite{Chan93} and $\sigma$ is as above.
\end{remark}

\begin{remark}
It is shown below in Lemma \ref{l:hstarint} that if
$H(x,y)=h_x(\abs{y})$ and $h_x\in\setH(\sigma)$ uniformly in $x\in\Rn$,
then $H$ satisfies \eqref{e:hrect}. For homogeneous functions the
situation is more complicated:

An interesting special case is $H(x,y)=\abs{\Omega(y)}$, $\Omega$ is
homogeneous of degree zero, and the set $S=S_\Omega$ is the star-shaped set
associated with $\Omega$ as in definition~\ref{d:S}. If the cover
$\{R_j\}$ is constructed as in remark \ref{r:coverconstr} (i.e., the
rectangles $\{R_j\}$ are the the same as $\{R_{m,k}\}$ after a
renumbering of the latter), condition \eqref{e:hrect} implies
$\Omega\in\Li(\Sn)$, as we now show.

Using the notation of remark~\ref{r:coverconstr} let $D_{m,k}$ be the
intersection of a certain cone inside $R_{m,k}$ with $\Sn$. These disks
were chosen to cover the set $\Theta_m$.  But then there exists a
constant $\beta\in(0,1)$ such that for each $m$ at least one of the
disks covering $\Theta_m$, say $D_{m,k_m}$,  has the property that
$\rho(\theta)>2^{m-1}$ on a subset of $D_{m,k_m}$ with surface measure
at least $\beta\abs{D_{m,k_m}}$ (we use $\abs{\cdot}$ also for surface
measure).   
Then, using $\abs{R_{m,k_m}} \approx 2^{nm}\abs{D_{m,k_m}}$, we have
\[
    \intave{R_{m,k_m}}{\abs{\Omega(y)} \,dy}
        \geq
    {c_n} \intave{D_{m,k_m}}{\abs{\Omega(\theta)}\,d\theta}.
\]
The right hand side is larger than $c_{n,\beta} 2^{n(m-1)}$, which
shows \eqref{e:hrect} is impossible when the set $S_\Omega$ is
essentially unbounded. Hence $\Omega\in \Li(\Sn)$.

If the covering is less efficient there are examples of
unbounded $\Omega$ that satisfy \eqref{e:hrect} and yield a bounded
operator $M_{S_\Omega,\Omega}$. In two dimensions let
\[
    \Omega(\theta)=\sum_{k=1}^\infty 2^{2k} \kar_{I_k}(\theta), 
	\quad \text{where $I_k = [2^{-3k}-2^{-5k},2^{-3k}]$,}
\]
and $R_j = [-2^j,2^j]\times[-2^{-2j},2^{-2j}]$, so that
$\sum_j\abs{R_j}<\infty$ and $\bigcup_jR_j\supset S_\Omega$.
A straightforward computation shows~$\Omega$ satisfies \eqref{e:hrect}.
Examples of weights that satisfy \eqref{e:AprRj} are $w(x)=\abs{x}^\alpha$ for
$-1<\alpha<p$. As a result $M_{S_\Omega,\Omega}$ is bounded on $\Lpw$.

It is of course easy to find examples of star-shaped sets $S$ and unrelated
homogeneous functions $H(x,y)=\abs{\Phi_x(y)}$, $\Phi(y)\in
L^\sigma(\Sn)$, that satisfy \eqref{e:hrect}: one idea is simply to
require $\Phi$ to be (uniformly) bounded along any unbounded arms of $S$.
\end{remark}

\begin{proof}[Proof of Theorem~\ref{t:MSh}]
The proof is almost identical with the one given in \cite{Chan93} in the
case $H\equiv1$, but since there are many differences in the details
we present the full argument:

Since the sets $R_j$ cover $S$ we have
$\kar_S\leq \sum_j \kar_{R_j}$, hence $M_{S,H}f \leq \sum_j M_{R_j,H}f$.
Let $\Lambda_j\colon\Rn\rightarrow\Rn$ be an invertible linear
transformation such that $R_j=\Lambda_j Q_1$, where~$Q_1$ is the unit
cube centered at the origin.
A change of coordinates gives
\[
    \int_{R_j} H(x,y) f(x-y)\,dy
    =
    \abs{\det\Lambda_j}
        \int_{Q_1} H(x,\Lambda_j y) f(\Lambda_j(\Lambda_j^{-1}x-y))
            \,dy,
\]
and letting
    $\Lambda_j f(x) = \abs{\det\Lambda_j} f(\Lambda_j x)$,
    $H_j(x,y) = H(x,\Lambda_j y)$,
we have 
\[
    M_{R_j,H}f(x) 
        \leq C \abs{\det\Lambda_j} (\Lambda_j^{-1} M_{H_j} \Lambda_j f)(x).
\]
Using $\snorm{\Lambda_j^{-1} g}{p,w} = \abs{\det\Lambda_j}^{-1}
\snorm{g}{p,\Lambda_j w}$ we get
\begin{align*}
    \snorm{M_{R_j,H}f}{p,w}
    &\leq
    C\abs{\det\Lambda_j}\, \snorm{\Lambda_j^{-1} M_{H_j} \Lambda_j f}{p,w}
    =
    C\snorm{M_{H_j} \Lambda_j f}{p,\Lambda_j w}   \\
    &\leq
    C\snorm{M_{H_j}}{L^p_{\Lambda_j w},L^p_{\Lambda_j w}}
    \snorm{\Lambda_j f}{p,\Lambda_j w} \\
    &=
    C\abs{\det\Lambda_j}\, \snorm{M_{H_j}}{L^p_{\Lambda_j w},L^p_{\Lambda_j w}}
        \snorm{f}{p,w}.
\end{align*}
Thus
$
    \snorm{M_{S,H}f}{p,w}
    \leq
    C\snorm{f}{p,w} \sum_j
        \abs{R_j} \,\snorm{M_{H_j}}{L^p_{\Lambda_j w},L^p_{\Lambda_j w}}
$,
since $\abs{R_j} = \abs{\det\Lambda_j}$.

The functions $H_j$ satisfy condition \eqref{e:hcube} uniformly in $j$:
Since for any $r>0$ there is $t>0$ such that
$tR_j \subset\Lambda_j B(0,r)\subset 2tR_j$, we get from \eqref{e:hrect} that
\begin{align*}
    \int_{\abs{y}<r} H_j(x,y)^\sigma\,dy
        &= \abs{\det\Lambda_j}^{-1}\int_{\Lambda_j B(0,r)} H(x,y)^\sigma\,dy\\
        &\leq C_H\abs{\det\Lambda_j}^{-1}\abs{\Lambda_j B(0,r)},
\end{align*}
which is equal to $C r^n$.
By a similar change of coordinates condition~\eqref{e:AprRj} is equivalent with 
\[
    \left(\intave{Q}{\Lambda_j w}\right)^{1/p}
    \left(\intave{Q}{(\Lambda_j w)^{-rp'/p}}\right)^{1/rp'} \leq \frac{K_j}{\abs{R_j}},
\]
for all cubes $Q\subset\Rn$.
By Theorem \ref{t:maxApr},
$\snorm{M_{H_j}}{L^p_{\Lambda_j w},L^p_{\Lambda_j w}} \leq C K_j/\abs{R_j}$.
Therefore $\snorm{M_{S,H}f}{p,w} \leq \bigr(C \sum_j K_j \bigl)\snorm{f}{p,w}$.
\end{proof}


\subsection{Vector-valued inequalities}

In this section we prove a generalization of the weighted vector-valued
inequality for the maximal function:

\begin{theorem}[Andersen and John \cite{Ande80}]\mylabel{t:AndJohn}
Let $1<p, q<\infty$ and suppose $w\in A_p$. There is a
constant $C_{p,q,w}$ such that
\[
    \Bigl\lVert\Bigl(\sum_j\abs{Mf_j}^q\Bigr)^{1/q}\Bigr\rVert_{p,w}
	\leq
    C_{p,q,w} \Bigl\lVert\Bigl(\sum_j\abs{f_j}^q\Bigr)^{1/q}\Bigr\rVert_{p,w},
\]
where $C_{p,q,w}$ depends on $w$ only through its $A_p$ constant.
\end{theorem}

The corresponding result for the operator $M_H$ is as follows:

\begin{theorem}\mylabel{t:maxvect}
Let $1<p, q<\infty$ and $w\in A_p$. There exists $\sigma>1$
such that if $H$ satisfies \eqref{e:hcube}
the vector-valued inequality
\begin{equation}\mylabel{e:maxvect}
    \Bigl\lVert\Bigl(\sum_j\abs{M_H f_j}^q\Bigr)^{1/q}\Bigr\rVert_{p,w}
	\leq
    C \Bigl\lVert\Bigl(\sum_j\abs{f_j}^q\Bigr)^{1/q}\Bigr\rVert_{p,w}
\end{equation}
holds for any sequence $\{f_j\}$.
\end{theorem}

\begin{proof}
By H\"older's inequality we have
\[
    \Bigl\lVert\Bigl(\sum_j\abs{M_H f_j}^q\Bigr)^{1/q}\Bigr\rVert_{p,w}
	\leq
    C_H^{1/\sigma}
    \Bigl\lVert\Bigl(\sum_j[M(\abs{f_j}^{\sigma'})]^{q/\sigma'}\Bigr)^{1/q}\Bigr\rVert_{p,w}
\]
and by Theorem~\ref{t:AndJohn} this is bounded by 
$C \bigl\lVert\bigl(\sum_j\abs{f_j}^q\bigr)^{1/q}\bigr\rVert_{p,w}$,
provided that $q/\sigma'>1$ and $w\in A_{p/\sigma'}$, both of which hold
when $\sigma$ is large enough (see~\cite{Coif74b,Muck72}).
\end{proof}
 

\section{Proof of the representation formula}\mylabel{s:srep}

Recall the class of functions $\setH(\sigma)$ from definition~\ref{d:H}.

\begin{lemma}\mylabel{l:Hlog}
Let $1\leq\sigma<\infty$, $h\in\setH(\sigma)$, and $C_h$ be as in the
definition of $\setH(\sigma)$. If $0<a<b<\infty$, then
$
    \int_a^b \abs{h(r)}^\sigma \,dr/r
	\leq C_h\lceil \log_2 b/a \rceil
$,
where $\lceil x\rceil $ is the ceiling of $x\in\mathbb{R}$, i.e., the
smallest integer greater than or equal to~$x$.
\end{lemma} 

\begin{proof}
Let $N = \lceil \log_2 b/a \rceil$, then
\[
    \int_a^b \abs{h(r)}^\sigma\, \frac{dr}r \leq
    \sum_{k=0}^{N-1} \int_{2^k a}^{2^{k+1}a}\abs{h(r)}^\sigma\, \frac{dr}r  \leq
    C_h N. \qed
\]
\noqed
\end{proof}

\begin{lemma}\mylabel{l:srep}
Let $\Omega\in L^1(\Sn)$ be positively homogeneous of degree $0$. Then
for all $y\in\Rn$
\begin{equation}\mylabel{e:srep:l}
    \int_0^\infty \frac1{t^n} \kar_{tS\setminus B(0,\epsilon)}(y) \frac{dt}t =
    \frac1n \kar_\compl{B(0,\epsilon)}(y) \frac{\abs{\Omega(y)}}{\abs{y}^n}
\end{equation}
in the sense that either both sides are finite and equal or they are both
infinite.
\end{lemma}

\begin{proof}
The identity follows from the fact that $y\in tS$ if and only if
$\abs{y}<t\abs{\Omega(y)}^{1/n}$.
\end{proof}

\begin{proof}[Proof of Theorem \ref{t:srep}]
Fix $x\in\Rn$ and $\epsilon>0$.
%
%
By Lemma \ref{l:srep} the operator $T_\epsilon$ can be written as
\begin{equation}\mylabel{e:srepB}
    T_\epsilon f(x) =
    n \int_\Rn \int_0^\infty \frac1{t^n} 
    \kar_{tS\setminus B(0,\epsilon)}(y)\frac{dt}t
    \sgn\Omega(y) h(\abs{y}) f(x-y) \,dy.
\end{equation}
Changing the order of integration---justification is given below---we obtain
\[
    T_\epsilon f(x) 
	= 
    n \int_0^\infty \frac1{t^n} \int_{tS\setminus B(0,\epsilon)} 
    \sgn\Omega(y) h(\abs{y}) f(x-y) \,dy \,\frac{dt}t,
\]
as was to be shown.

To justify the change in the order of integration
we need to show that
\begin{multline}\mylabel{e:srepA}
    n \int_{\Rn} \int_0^\infty \frac1{t^n} \kar_{tS\setminus B(0,\epsilon)}(y) 
    \frac{dt}t \abs{h(\abs{y})f(x-y)}\,dy \\
    =
    \int_{\abs{y}>\epsilon}\frac{\abs{\Omega(y)}}{\abs{y}^n}\abs{h(\abs{y})f(x-y)}\,dy
\end{multline}
is finite.
%
%
A change into polar coordinates $y=r\theta$ with $r=\abs{y}$ and 
$\theta=y/\abs{y}\in\Sn$ shows the right-hand side of \eqref{e:srepA} is
\begin{equation}\mylabel{e:srepD}
    \int_\Sn \abs{\Omega(\theta)} \int_\epsilon^\infty
	\frac{\abs{h(r)}}r
        \abs{f(x-r\theta)} \,dr \,d\theta,
\end{equation}
where we used $d\theta$ for the surface measure on $\Sn$.
We estimate \eqref{e:srepD} by
\[
    \snorm{f}{\infty}
    \int_\Sn \abs{\Omega(\theta)} \,d\theta 
    \int_\epsilon^{R+\abs{x}} \abs{h(r)}\,\frac{dr}{r},
\]
where $R$ is such that $\supp f\subset B(0,R)$. By Lemma \ref{l:Hlog}
the second integral is finite for all $x\in\Rn$, thus also \eqref{e:srepA}
is finite everywhere.

Since the case $f\in L^p(\Rn)$, $1<p<\infty$, is not used in the rest
of the paper, we only sketch its proof: We let $Q\subset\Rn$ be an
arbitrary cube and compute the $L^p(Q)$ norm of \eqref{e:srepD}. Using
Minkowski's inequality, $h\in\setH(p')$, and $f\in L^p(\Rn)$, we can
estimate~\eqref{e:srepD} by $C(Q,h,\epsilon,f) \int_\Sn
\abs{\Omega(\theta)}\,d\theta<\infty$. Since the cube $Q$ is
arbitrary, this shows that \eqref{e:srepA} is finite almost
everywhere.
\end{proof}

\section{Preliminaries}\mylabel{s:prelim}

\subsection{Weighted Littlewood-Paley theory}\mylabel{sec:LP}

For proofs of the following facts see \cite{Kurt79} and \cite{Wats97}.

Let $\Psi\in C_0^\infty((\frac12,2))$ be such that
$\Psi\geq0$ and
$\sum_{j\in\Z} \Psi^2(2^jt) = 1$ for $t>0$.
Let $\fourier{\psi}(\xi) = \Psi(\abs{\xi})$ and 
define $Q_j f = \psi_j * f$, where
$\psi_j(x) = 2^{-nj} \psi(2^{-j}x)$. 

Let $1<p<\infty$ and $w\in A_p$. Then there exists $C_{p,w}>0$ such that
the following estimates hold:
\begin{allowdisplaybreaks}
\begin{gather}
    \mylabel{e:LP1}
    C_{p,w}^{-1} \snorm{f}{p,w} \leq \bsqnorm{\sum_j \abs{Q_j f}^2}
        \leq C_{p,w} \snorm{f}{p,w}, \\
    \mylabel{e:LP1p}
    C_{p,w}^{-1} \snorm{f}{p,w} \leq \bsqnorm{\sum_j \abs{Q_j^2 f}^2}
        \leq C_{p,w} \snorm{f}{p,w},  \\
    \mylabel{e:LP2}
    \mnorm{\sum_j Q_j f_j}{p,w}   \leq C_{p,w} \bsqnorm{\sum_j \abs{Q_j f_j}^2},
	 \\
    \mylabel{e:LP2p}
    \mnorm{\sum_j Q_j^2 f_j}{p,w}   \leq C_{p,w} \bsqnorm{\sum_j \abs{Q_j^2 f_j}^2}.
\end{gather}
\end{allowdisplaybreaks}

The identity $\sum_j Q_j^2 = \text{Id}$ holds in many senses, e.g., for smooth
functions $f\in L^2$
\begin{equation}\mylabel{e:LPid}
    \lim\limits_{k,l\rightarrow\infty} \sum_{j=-k}^l Q_j^2 f(x)
                = f(x)
\end{equation}
for all $x\in\Rn$.

\subsection{Lemmata on homogeneous and radial functions}

In the following lemmata $\Omega$ is always positively homogeneous of
degree~$0$ and $S=S_\Omega$ is the star-shaped set associated with
$\Omega$ as in definition \ref{d:S}.

\begin{allowdisplaybreaks}
\begin{lemma}\mylabel{l:Sints}
The following equalities and estimates hold:
\begin{align}
    \mylabel{e:Sint}
    \abs{S} = \int_S \, dy &= \frac1n 
                \int_\Sn \abs{\Omega(\theta)}\,d\theta,
                \\
    \mylabel{e:Ssgnint}
    \int_S \sgn\Omega(y)\, dy &= \frac1n \int_\Sn \Omega(\theta)\,d\theta,
                \\
    \mylabel{e:Slogpint}
    \int_S \logp \abs{y} \, dy &\leq \frac1{n^2} \snorm{\Omega}{\LlogL},  	
                \\
    \mylabel{e:Slogint}
    \int_S \bigl|\log \abs{y}\bigr| \, dy 
        &\leq \frac1{n^2} (\snorm{\Omega}{\LlogL} + 1),
                \\
    \mylabel{e:mSmsum}
    \sum_{m=0}^\infty (m+1)\abs{S_m} &\leq c_n \snorm{\Omega}{\LlogL},
\end{align}    
for some constant $c_n$ depending only on the dimension~$n$.
\end{lemma}
\end{allowdisplaybreaks}

\begin{proof}
Using polar coordinates $y=r\theta$ with $r>0$, 
$\theta\in\Sn$,
and recalling that $\rho(\theta) = \abs{\Omega(\theta)}^{1/n}$, we have
\[
    \int_S \, dy = \int_\Sn \int_0^{\rho(\theta)} r^{n-1}\,dr\,d\theta
        = \frac1n \int_\Sn \abs{\Omega(\theta)}\,d\theta.
\]
This proves \eqref{e:Sint}. The proofs of \eqref{e:Ssgnint} through
\eqref{e:Slogint} are similar.

To prove estimate \eqref{e:mSmsum} recall the sets $\Theta_m$ from
\eqref{e:Thmdef} and notice that $\abs{S_m} = \frac1n\int_{\Theta_m}
\abs{\Omega(\theta)}\,d\theta$ by \eqref{e:Sint}. This gives
\[
\sum_{m=0}^\infty (m+1)\abs{S_m} 
    = \frac1n \int_{\Theta_0} \abs{\Omega(\theta)}\,d\theta + \frac1n 
       \sum_{m=1}^\infty \int_{\Theta_m} (m+1)\abs{\Omega(\theta)} \,d\theta.
\]
When $m\geq1$, $\rho(\theta)>2^{m-1}\geq1$ on $\Theta_m$;
thus $m< 1+\logpb \abs{\Omega(\theta)}^{1/n}$
and so
\[
\sum_{m=0}^\infty (m+1)\abs{S_m} 
    \leq \frac1n \int_\Sn
        \abs{\Omega(\theta)}\,d\theta + \frac1n \int_\Sn
        \abs{\Omega(\theta)} (2+\logpb \abs{\Omega(\theta)}^{1/n}) \,d\theta,
\]
which is bounded by $c_n \snorm{\Omega}{\LlogL}$.
\end{proof}

We will use the following characterizations of the class $\setH(\sigma)$
interchangeably:

\begin{lemma}\mylabel{l:Heq}
For fixed $1\leq\sigma<\infty$ the following are equivalent:
\begin{enumerate}
\item[{\upshape(a)}]
    there exists $C>0$ such that $\int_R^{2R} \abs{h(r)}^\sigma\,dr \leq CR$
    for all $R>0$.
\item[{\upshape(b)}]
    there exists $C>0$ such that $\int_0^R \abs{h(r)}^\sigma\,dr \leq CR$
    for all $R>0$.
\item[{\upshape(c)}]
    there exists $C>0$ such that $\int_R^{2R} \abs{h(r)}^\sigma\,\frac{dr}r
	 \leq C$
    for all $R>0$.
\end{enumerate}
\end{lemma}

The only non-trivial implication is from (a) to (b), which follows by
writing the integral over $[0,R]$ as the sum of integrals over
$[2^{-k-1}R,2^{-k}R]$, $0\leq k<\infty$.

\begin{lemma}\mylabel{l:hSint}
Suppose that $h\in\setH(1)$
vanishes in a neighborhood of the origin. Then there are
constants $C_1(h)>0$ and $C_2(h)>0$ such that
\begin{equation}\mylabel{e:hint}
    \int_0^R \lvert h(r)\rvert \frac{dr}r \leq C_1(h) + C_2(h) \logp R,
\end{equation}
and a constant $C(h,n)>0$ such that
\begin{equation}\mylabel{e:hSint}
    \int_0^1 \int_S \lvert h(t\abs{y}) \rvert \, dy\frac{dt}t 
    \leq
    C(h,n) \snorm{\Omega}{\LlogL}.
\end{equation}
\end{lemma}

\begin{proof}[Sketch of Proof]
The first part of the lemma follows immediately from Lemma~\ref{l:Hlog}.
For the second part change the order of integration in~\eqref{e:hSint},
change variables in the $t$-integral, and use~\eqref{e:hint},
\eqref{e:Sint}, and~\eqref{e:Slogpint}
\end{proof}

\begin{lemma}\mylabel{l:hstarint}
Let $1\leq\sigma<\infty$, $h\in\setH(\sigma)$,
and $E\subset\Rn$ be any measurable set star-shaped about
the origin. Then
\[
    \int_E \abs{h(t\abs{y})}^\sigma\,dy \leq c_nC_h \abs{E},
        \quad \text{for any $t>0$,}
\]
where $C_h$ is as in the definition of $h\in\setH(\sigma)$
\textup{(}definition~\textup{\ref{d:H})}.
\end{lemma}

\begin{proof}
In polar coordinates
\[
    \int_E \abs{h(t\abs{y})}^\sigma\,dy
    = \int_\Sn \int_0^{\rho(\theta)} \abs{h(tr)}^\sigma r^{n-1}\,dr\,d\theta,
\]
where $\rho(\theta)$ is the boundary function of $E$, that is, up to a
set of measure zero $E$ is given by $\{x\in\Rn\colon \abs{x}\leq\rho(x)\}$.
Now the $r$-integral equals
\[
    \sum_{k=0}^\infty \int_{2^{-k-1}\rho(\theta)}^{2^{-k}\rho(\theta)}
         \abs{h(tr)}^\sigma r^{n-1}\,dr 
    \leq \sum_{k=0}^\infty 2^{n(-k-1)}\rho(\theta)^n
        \int_{2^{-k-1}\rho(\theta)}^{2^{-k}\rho(\theta)} \abs{h(tr)}^\sigma
        \,\frac{dr}r.
\]
Lemma \ref{l:Heq} and a change of coordinates show that the last
integral is bounded by a constant,
%
%
%
%
thus
$\int_0^{\rho(\theta)} \abs{h(tr)}^\sigma r^{n-1}\,dr \leq
c_nC_h \rho(\theta)^n$.
Hence
\[
    \int_E \abs{h(t\abs{y})}^\sigma\,dy \leq c_nC_h \int_\Sn\rho(\theta)^n\,d\theta
        = c_nC_h \abs{E}. \qed
\]
\noqed
\end{proof}


\section{Proof of the weighted norm estimates}\mylabel{s:hbded}

We prove Theorems \ref{t:Thbded} and \ref{t:Thinfbded} in this section. 
We first show that the uniform boundedness of the truncated operators
$T_\epsilon$ implies that $T$ is bounded: 
let $f,g\in C_0^1(\Rn)$ and use the uniform boundedness
to get for any $\epsilon>0$,
\[
    \abs{\ip{Tf}{g}} 
    \leq \abs{\ip{(T-T_\epsilon)f}{g}} + C \snorm{f}{\Lpw}\snorm{g}{(\Lpw)'},
\]
where $C=\sup_{\epsilon>0} \snorm{T_\epsilon}{\Lpw,\Lpw}$. Since by
remark~\ref{r:limTheps} $T_\epsilon f$ converges to $Tf$ uniformly, we
have $\abs{\ip{(T-T_\epsilon)f}{g}} \rightarrow 0$ as $\epsilon\tozerop$.
This shows $T$ is bounded on $\Lpw$ with operator norm at most $C$.

We will prove Theorem \ref{t:Thbded} only for $1<p\leq 2$. The
case $p>2$ follows by a standard duality argument: the adjoint
$T^*$ is essentially of the same form as $T$ and if the weight $w$
satisfies \eqref{e:ThcondCb} then the dual weight $w^{-p'/p}$ satisfies
\eqref{e:ThcondCa} with $p'$ in place of $p$. Hence $T^*$ is bounded on
$L^{p'}_{w^{-p'/p}}$ and so $T$ is bounded on $\Lpw$.

To simplify notation we leave out $\epsilon$ and instead assume
that $h$ vanishes in some neighborhood of the
origin. This corresponds to truncated operators by replacing an
arbitrary $h(r)$ by $h_\epsilon(r) = h(r)\kar_{(\epsilon,\infty)}(r)$.
We show that the operator norm depends on $h$ only through the
constants $C_h$ in the definition of $h\in\setH(\sigma)$. In particular,
the norm is independent of the support of $h$. Since
$\abs{h_\epsilon(r)}\leq \abs{h(r)}$, $r>0$, the functions $h_\epsilon$
belong to $\setH(\sigma)$ with the same constants as $h$, and therefore the
proof shows the truncated operators are uniformly bounded on $\Lpw$.

We concentrate on the proof of Theorem~\ref{t:Thbded} and indicate at
the end of this section how the argument is modified to prove
Theorem~\ref{t:Thinfbded}.  From now on, throughout the proof of
Theorem~\ref{t:Thbded}, we assume $h$ vanishes in a neighborhood of
the origin.

\subsection{The decompositions}

Corresponding to the decomposition of $S$ into the sets $S_m$,
see~\eqref{e:Smdef}, we define the averaging operators
\begin{equation}\mylabel{e:atmdef}
    A_t^mf(x) = a_t^m * f(x), \quad\text{where}\quad
    a_t^m(x) = \frac1{t^n} \kar_{tS_m}(x) h(\abs{x}) \sgn\Omega(x),
\end{equation}
and we also use
\begin{equation}\mylabel{e:bsmdef}
    B_s^mf(x) = b_s^m * f(x), \qquad\text{where}\quad
    b_s^m(x) = \int_{2^{s-1}}^{2^s} a_t^m(x) \frac{dt}t.
\end{equation}
Since $h$ vanishes near the origin the representation formula
\eqref{e:srepY} gives
$
    Tf(x) = n \int_0^\infty A_t f(x)\,dt/t$,
where, see \eqref{e:srepX},
$
    A_t f(x) = a_t * f(x)$ and
$    a_t(x) = t^{-n} \kar_{tS}(x) h(\abs{x}) \sgn\Omega(x)$.
For test functions $f\in C_0^1(\Rn)$ we write
\begin{equation}\mylabel{e:decomTsum}
    Tf(x) = n \sum_{s\in\Z,m\geq0} B_s^m f(x).
\end{equation}
This involves a change in the order of integration and summation,
which is justified by the following lemma. The lemma also shows the
order of summation does not matter.

\begin{lemma}\mylabel{l:Bsmbded}
If $\Omega$ is as in Theorem \textup{\ref{t:Thbded}}, $h\in\setH(2)$,
and $f\in C_0^1(\Rn)$, then
\[
    K(x)=
    \sum_{s\in\Z,m\geq0}
        \int_{2^{s-1}}^{2^s}
            \frac1{t^n} \int_{tS_m} \mabs{f(x-y)h(\abs{y})}\,dy\,\frac{dt}t
\]
is bounded on $\Rn$.
\end{lemma}

\begin{proof}
Evaluating the double sum gives
\[
    K(x)=
    \int_0^\infty 
        \frac1{t^n} \int_{tS}
        \abs{f(x-y) h(\abs{y})}\,dy \frac{dt}t.
\]
Now
\[
    \int_0^1 \frac1{t^n} \int_{tS}
        \abs{f(x-y) h(\abs{y})} \,dy \frac{dt}t
    \leq
    \snorm{f}{\infty} \int_0^1 \int_S \abs{h(t\abs{y})} \,dy\frac{dt}t,
\]
and by Lemma \ref{l:hSint} this is finite. On the other hand
\begin{multline*}
    \int_1^\infty \frac1{t^n} \int_{tS} \abs{f(x-y) h(\abs{y})} \,dy \frac{dt}t\\
    \leq 
    \left( \int_1^\infty \int_{tS} \abs{f(x-y)}^2 \,dy \frac{dt}{t^2} \right)^{1/2}
    \left( \int_1^\infty \int_{tS} \abs{h(\abs{y})}^2 \,dy \frac{dt}{t^{2n}} \right)^{1/2}.
\end{multline*}
The first factor is bounded by
$
    \left( \int_1^\infty \int_{\Rn} \abs{f(x-y)}^2 \,dy \frac{dt}{t^2} \right)^{1/2}
    = \snorm{f}{2}.
$
Using Lemma \ref{l:hstarint}, $\int_{tS} \abs{h(\abs{y})}^2 \,dy \leq C\abs{tS}$,
when $h\in\setH(2)$.
Therefore the second factor is at most
$    C \left( \int_1^\infty \abs{tS} dt/t^{2n} \right)^{1/2}
     \leq C_n \abs{S} < \infty$.
\end{proof}

When $f\in C_0^1$, $\int_0^\infty A_t f(x)\frac{dt}t =
\sum_{s,m} B_s^m f(x)$ is a smooth $L^2$ function (recall that $h$
vanishes near $0$). Let $Q_j$ be a Littlewood-Paley operator as in
Section~\ref{sec:LP}. By \eqref{e:LPid}, $\sum_{s,m} B_s^m f = \sum_j
Q_j^2 \sum_{s,m} B_s^m f$ pointwise and by the Lebesgue dominated
convergence theorem and Lemma~\ref{l:Bsmbded} we have
$
    Q_j^2 \sum_{s,m} B_s^m f =
        \sum_{s,m} Q_j^2 B_s^m f.
$ 
Hence
\begin{align*}
    \int_0^\infty A_t f(x)\frac{dt}t
        &= \sum_j Q_j^2 \sum_{s\in\Z,m\geq0} B_{s+j-m}^m f(x) \\
        &= \sum_j \sum_{s\in\Z,m\geq0} Q_j^2 B_{s+j-m}^m f(x).
\end{align*}
Therefore we may use the decomposition
$
    \int_0^\infty A_t f(x) \,dt/t = \I + \II + \III$,
where
\begin{allowdisplaybreaks}
\begin{align*}
   \I &= \sum_{j\in\Z} \sum_{m\geq0} \sum_{0\leq s\leq Nm} Q_j^2
        B_{s+j-m}^m f, \\
  \II &= \sum_{j\in\Z} \sum_{m\geq0} \sum_{s>Nm} Q_j^2
        B_{s+j-m}^m f, \\
 \III &= \sum_{j\in\Z} \sum_{m\geq0} \sum_{s<0} Q_j^2
        B_{s+j-m}^m f,
\end{align*}
for some $N\in\Z_+$ to be chosen later. 
\end{allowdisplaybreaks}

It is enough to show that each of the terms $\I$, $\II$, and $\III$
defines an operator bounded in $\Lpw$ norm. The estimate for $\I$
involves condition \eqref{e:ThcondC} on the weight. On the other hand
we show $\II$ and $\III$ are bounded on $\Lpw$ for any $A_p$ weight~$w$.
The proof for $\II$ also fixes the value of $N$.

\subsection{Term $\I$}

Let
\[
    G_j f(x) = \sum_{m=0}^\infty \int_{2^{j-m-1}}^{2^{j+(N-1)m}}
        A_t^m f(x) \frac{dt}t,
\]
where $A_t^m$ is defined in \eqref{e:atmdef}. Then $\I=\sum_j Q_j^2 G_j f$.
We will prove a square norm inequality
\begin{equation}\mylabel{e:Gjsqnorm}
    \bsqnorm{ \sum_j \abs{G_j f_j}^2 } 
                \leq C_{p,w,N} \bsqnorm{ \sum_j \abs{f_j}^2 }
\end{equation}
for all $N\in\Z_+$ and $1<p\leq2$,
where $C_{p,w,N}$ is independent of the sequence
$\{f_j\}$. Then, since $w\in A_p$ by remark \ref{r:wAp},
the Littlewood-Paley estimate \eqref{e:LP2p} gives
\[
    \snorm{\I}{p,w} =
    \bnorm{\sum_j Q_j^2 G_j f}{p,w} \leq
    C \bsqnorm{ \sum_j \abs{Q_j^2 G_j f} ^2 }.
\]
Since both $G_j$ and $Q_j$ are convolution operators, they commute. Hence
\[
    \snorm{\I}{p,w} \leq
    C \bsqnorm{ \sum_j \abs{ G_j Q_j^2 f }^2 } \leq
    C \bsqnorm{ \sum_j \abs{ Q_j^2 f }^2 },
\]
where the second inequality follows from \eqref{e:Gjsqnorm}. Again by
Littlewood-Paley theory, equation \eqref{e:LP1p} this time, we get
$\snorm{\I}{p,w} \leq C \snorm{f}{p,w}$. 

The inequality \eqref{e:Gjsqnorm} is shown by proving the following two
estimates: 
\begin{equation}\mylabel{e:Gjpw}
    \snorm{G_j f}{p,w} \leq C \snorm{f}{p,w},
\end{equation}
and
\begin{equation}\mylabel{e:Gjsuppw}
    \mnorm{\sup_j \abs{G_j f_j}}{p,w} \leq C \mnorm{\sup_j \abs{f_j}}{p,w},
\end{equation}
where $C$ is a constant independent of $f$ and $\{f_j\}$. Then
\eqref{e:Gjpw} implies the vector-valued inequality
\begin{equation}\mylabel{e:Gjppw}
     \biggl\lVert \biggl( \sum_j \abs{G_j f_j}^p \biggr)^{1/p} \biggr\rVert_{p,w}
        \leq
     C \biggl\lVert \biggl( \sum_j \abs{f_j}^p \biggr)^{1/p} \biggr\rVert_{p,w}.
\end{equation}
Interpolation in the sequence norm
between \eqref{e:Gjsuppw} and \eqref{e:Gjppw} gives
\eqref{e:Gjsqnorm} for $1<p\leq 2$.

We begin proving \eqref{e:Gjpw} by first defining the positive operator
\[
    A_{+t}^{m,k} f(x) = \frac1{t^n} \int_{tR_{m,k}}
                                \abs{h(\abs{y})} f(x-y)\,dy,
\]
where $\{R_{m,k}\}$ is the stratified starlike cover of $S$.
Since $S_m \subset \bigcup_k R_{m,k}$ we have
$\abs{A_t^m f(x)} \leq \sum_k A_{+t}^{m,k}(\abs{f})(x)$ and
by applying Minkowski's inequalities we get
\begin{equation}\mylabel{e:GjpwA}
    \snorm{G_j f}{p,w} \leq
        \sum_{m\geq0} \int_{2^{j-m-1}}^{2^{j+(N-1)m}} \sum_k
                \snorm{A_{+t}^{m,k}(\abs{f})}{p,w} \,\frac{dt}t.
\end{equation}
H\"older's inequality
and a change of coordinates gives
\newcommand{\xtRmk}[1][x]{R_{#1}}
\begin{multline*}
    A_{+t}^{m,k}(\abs{f})(x) =
        \frac1{t^n} \int_{\xtRmk}
            \abs{h(\abs{x-y})} \abs{f(y)} w(y)^{1/p} w(y)^{-1/p}\,dy
    \\ \leq
        \frac1{t^n}
        \left( \int_{\xtRmk} \abs{f(y)}^{p} w(y)\,dy \right)^{1/p}
        \left( \int_{\xtRmk} \abs{h(\abs{x-y})}^{p'} w(y)^{-p'/p}\,dy \right)^{1/p'},
\end{multline*}
where we used $\xtRmk=x-tR_{m,k}$ for fixed $t$, $m$, and $k$ to simplify
notation. Another application of H\"older's inequality gives
\begin{multline}\mylabel{e:anotherholder}
    \left( \int_{\xtRmk} \abs{h(\abs{x-y})}^{p'} w(y)^{-p'/p}\,dy \right)^{1/p'}
            \\ \leq
    \left( \int_{tR_{m,k}} \abs{h(\abs{y})}^{r'p'}\,dy \right)^{1/r'p'}
    \left( \int_{\xtRmk} w(y)^{-rp'/p}\,dy \right)^{1/rp'},
\end{multline}
and when $\sigma\geq r'p'$ Lemma \ref{l:hstarint} shows this is bounded by
\[
    C \abs{tR_{m,k}}^{1/p'}
        \left( \intave{\xtRmk}{w(y)^{-rp'/p}\,dy} \right)^{1/rp'}.
\]
%
%
%
%
Thus $\snorm[p]{A_{+t}^{m,k}(\abs{f})}{p,w}$ is bounded by
\[
    C \frac{\abs{R_{m,k}}^{p/p'}}{t^n} \int_\Rn
        \int_{\xtRmk} \abs{f(y)}^{p} w(y)\,dy
        \left( \intave{\xtRmk}{w(z)^{-rp'/p}\,dz} \right)^{p/rp'}
    w(x) \, dx,
\]
and a change in the order of integration shows this equals
\[
    C \frac{\abs{R_{m,k}}^{p/p'}}{t^n}
        \int_\Rn \biggl[ \abs{f(y)}^{p} w(y)
        \int_{\xtRmk[y]} w(x)
        \left( \intave{\xtRmk}{w(z)^{-rp'/p}\,dz} \right)^{p/rp'}
    \, dx \biggr] \,dy,
\]
where we also used the symmetry of the rectangle $R_{m,k}$.
Now $z\in \xtRmk=x-tR_{m,k}$ and $x\in \xtRmk[y]=y-tR_{m,k}$ imply
$z\in y-2tR_{m,k}$. Hence we get the bound
\begin{multline*}
    C \abs{R_{m,k}}^{p}
        \int_\Rn \biggl[ \abs{f(y)}^{p} w(y)  \;
        \frac1{\abs{R_y}} \int_{ y-2tR_{m,k} }{ w(x)\,dx } \\
        \times
        \biggl( \frac1{\abs{R_y}}
            \int_{ y-2tR_{m,k} }{ w(z)^{-rp'/p}\,dz } \biggr)^{p/rp'}
        \biggr]\, dy.
\end{multline*}
By condition \eqref{e:ThcondCa} on the weight $w$ this is bounded by
$C \bigl(K_{m,k} \snorm{f}{p,w}\bigr)^p$,
for $1<p\leq 2$. 
Thus $\snorm{A_{+t}^{m,k}(\abs{f})}{p,w} \leq C K_{m,k} \snorm{f}{p,w}$.
Substituting into \eqref{e:GjpwA} gives
\[
    \snorm{G_j f}{p,w} \leq
        C \sum_{m\geq0} \sum_k (Nm+1) K_{m,k} \snorm{f}{p,w},
\]
which, by \eqref{e:Thsumcond}, is bounded by $C \snorm{f}{p,w}$. This
proves \eqref{e:Gjpw}.

To prove \eqref{e:Gjsuppw} let 
$
    Gf = \sum_{m\geq0} (m+1)M_{S_m,\abs{h}} (f)
$,
where, with slight abuse of notation, $M_{S_m,\abs h}$ is the maximal
operator of Theorem~\ref{t:MSh} with the function
$H(x,y)=\abs{h(\abs{y})}$ in the kernel.
Then
\begin{multline}\mylabel{e:GjsuppwA}
    \abs{G_j f(x)} \leq 
        \sum_{m\geq0} \int_{2^{j-m-1}}^{2^{j+(N-1)m}} \frac{dt}t  \,
            \sup_{t>0} \frac1{t^n} \int_{tS_m}
                \abs{h(\abs{y}) f(x-y)} \,dy
    \\ =
        C\sum_{m\geq0} (Nm+1)  M_{S_m,\abs{h}} (f)(x)
    \leq
        C_N G(f)(x).
\end{multline}
Since $\{R_{m,k}\}_{m,k}$ is a stratified starlike cover of the set $S$, for
fixed $m$ the collection $\{R_{m,k}\}_k$ is a starlike cover of $S_m$.
Then Lemma \ref{l:hstarint} shows
$\abs{h(\abs{\cdot})}$ satisfies \eqref{e:hrect} on the rectangles
$\{R_{m,k}\}$. 
Theorem \ref{t:MSh} gives $\snorm{M_{S_m,\abs{h}} (f)}{p,w}
\leq C \sum_k K_{m,k} \snorm{f}{p,w}$, thus
\begin{equation}\mylabel{e:GjsuppwB}
    \snorm{Gf}{p,w} \leq
        C \sum_{m\geq0} (m+1) \sum_k K_{m,k} \snorm{f}{p,w}
        \leq C \snorm{f}{p,w},
\end{equation}
where the last inequality follows from \eqref{e:Thsumcond}. We get from
\eqref{e:GjsuppwA} 
\[
    \mnorm{\sup_j \abs{G_j f_j}}{p,w} \leq
        C \mnorm{\sup_j Gf_j}{p,w} \leq
        C \mnorm{G(\sup_j \abs{f_j})}{p,w},
\]
and then \eqref{e:GjsuppwB} gives the bound $C \smnorm{\sup_j \abs{f_j}}{p,w}$,
which proves \eqref{e:Gjsuppw} and completes the proof of the
boundedness of the operator defined by term~$\I$.

\subsection{Term $\II$}

The norm estimate for $\II$ is based on first proving a good
unweighted $L^2$ estimate using Fourier transform techniques. Next
bounding the terms in $\II$ by maximal functions yields a crude
weighted estimate. The final estimate is obtained by
interpolating with change of measure.

To find a good unweighted estimate for the terms in $\II$ we begin by
estimating the Fourier transform of $a_t^m$. We write
$\fourier{a_t^m}(\xi)$ in the form
\[
    \fourier{a_t^m}(\xi) =
    \int_\Rn a_t^m(x) e^{-2\pi i x\cdot\xi} dx
    =
    \frac1{t^n} \int_{tS_m} h(\abs{x}) \sgn\Omega(x) e^{-2\pi i x\cdot\xi} dx.
\]
Using polar coordinates and making a change in the order of
integration this is equal to
\[
    \int_{\Theta_m} \int_0^{\rho(\theta)} h(tr) \sgn\Omega(\theta)
        e^{-2\pi i tr\theta\cdot\xi} r^{n-1} dr\,d\theta
                =
    \int_0^{2^m} h(tr) I_r(\xi) r^{n-1} dr,
\]
where, for fixed $t$ and $m$,
$
    I_r(\xi) = \int_{\Theta_m(r)}
        \sgn\Omega(\theta)
        e^{-2\pi i tr\theta\cdot\xi} \,d\theta
$
and $\Theta_m(r) = \{\theta\in\Theta_m\colon \rho(\theta)>r\}$. Hence
by Schwarz's inequality
\[
    \abs{\fourier{a_t^m}(\xi)}^2 \leq
    \int_0^{2^m} \abs{h(tr)}^2 r^{n-1} dr
    \int_0^{2^m} \abs{I_r(\xi)}^2 r^{n-1} dr.
\]
Note that $\int_0^{2^m} \abs{h(tr)}^2 r^{n-1} dr =
c_n \int_{B(0,2^m)} \abs{h(t\abs{y})}^2\,dy \leq c_n 2^{mn}$,
the inequality following from Lemma \ref{l:hstarint}.
To estimate the second factor we use some ideas from~\cite{Duoa86}
and write the square of $\abs{I_r}$ in the form
\[
    \abs{I_r(\xi)}^2 = \int_{\Theta_m(r)} \int_{\Theta_m(r)}
        \sgn\Omega(\theta) \overline{\sgn\Omega(\omega)}
            e^{-2\pi i tr(\theta-\omega)\cdot\xi}
        \,d\theta\,d\omega,
\]
so that $\int_0^{2^m} \abs{I_r(\xi)}^2 r^{n-1} dr$ is equal to
\begin{equation}\mylabel{e:IIvaliA}
    \int_{\Theta_m} \int_{\Theta_m}
        \sgn\Omega(\theta) \overline{\sgn\Omega(\omega)}
        \biggl[
            \int_0^{\bmin{\rho(\theta)}{\rho(\omega)}}
                e^{-2\pi i tr(\theta-\omega)\cdot\xi}
                    r^{n-1} dr
        \biggr]
        \,d\theta\,d\omega,
\end{equation}
where $\bmin{\rho(\theta)}{\rho(\omega)} =
\min\{\rho(\theta),\rho(\omega)\}$. 

A direct estimation combined with an integration by parts shows, for $a>0$, $b\not=0$,
$
    \left|\int_0^a e^{-ibr} r^{n-1} dr\right| 
	\leq
    C_\alpha {a^n}/{(a\abs{b})^\alpha},
$
where $0<\alpha<1$ and $C_\alpha>0$. Fix such an~$\alpha$. In particular
we get that
$
    \left|\int_0^a e^{-ibr} r^{n-1} dr\right| 
	\leq
    C {a_0^n}/{(a_0\abs{b})^\alpha},
$
when $a_0\leq a \leq 2a_0$, $a_0>0$.
To apply this to \eqref{e:IIvaliA} note that $\rho(\theta)$ and
$\rho(\omega)$ are between $2^{m-1}$ and~$2^m$ on~$\Theta_m$.
Taking $b=t(\theta-\omega)\cdot\xi$
the absolute value of the expression inside the brackets in \eqref{e:IIvaliA}
is bounded by
$
    C {2^{mn}}/{\abs{2^m t(\theta-\omega)\cdot\xi}^\alpha}.
$
Hence
\[
    \int_0^{2^m} \abs{I_r(\xi)}^2 r^{n-1} dr
    \leq
    C 2^{mn} \int_\Sn \int_\Sn
        \frac1{\abs{2^m t(\theta-\omega)\cdot\xi}^\alpha}
    \,d\theta\,d\omega.
\]
Since $0<\alpha<1$, $\abs{\theta-\omega}^{-\alpha}$ is integrable on
$\Sn\times\Sn$. Hence the above expression 
is bounded by $C 2^{mn} \abs{2^mt\xi}^{-\alpha}$.
Thus we have shown
$
    \abs{\fourier{a_t^m}(\xi)} \leq
    C {2^{mn}}/{\abs{2^mt\xi}^{\alpha/2}}
$,
which gives
\[
    \abs{\fourier{b_s^m}(\xi)} =
    \biggl| \int_{2^{s-1}}^{2^s} \fourier{a_t^m}(\xi) \frac{dt}t
    \biggr|
    \leq
    C \frac{2^{mn}}{\abs{2^{m+s}\xi}^{\alpha/2}} .
\]

Recall that $Q_j$ is convolution by $\psi_j(x) = 2^{-nj} \psi(2^{-j}x)$ and
the support of $\fourier\psi(\xi)$ is contained in $1/2\leq\abs{\xi}\leq2$,
see Section~\ref{sec:LP}.
Since $\snorm{\fourier{\psi}_j}{\infty} = \snorm{\fourier{\psi}}{\infty}$ 
and the support of $\fourier{\psi}_j$ is contained in the annulus
$\mathcal{A}_j = \{\xi\in\Rn\colon 2^{-j-1} \leq \abs{\xi} \leq 2^{-j+1} \}$,
\[
    \abs{(Q_j B_{s+j-m}^m f)\foursup{}(\xi)}
    =
    \abs{\fourier{\psi}_j(\xi) \fourier{b}_{s+j-m}^m(\xi) \fourier{f}(\xi)}
    \leq
        C \frac{2^{mn}}{\abs{2^{s+j}\xi}^{\alpha/2}}
        \abs{\fourier{f}(\xi)} \kar_{\mathcal{A}_j}(\xi).
\]
If $\xi\in\mathcal{A}_j$, then $\abs{2^{s+j}\xi}^{-\alpha/2}
\leq 2^{-s\alpha/2+\alpha/2}$, and thus
\[
    \abs{(Q_j B_{s+j-m}^m f)\foursup{}(\xi)}
    \leq
        C 2^{mn} 2^{-s\alpha/2} 
        \abs{\fourier{f}(\xi)}.
\]
Hence, by applying Plancherel's theorem twice,
we arrive at the inequality
\begin{equation}\mylabel{e:good2est}
    \bsqnorm[2]{\sum_j \abs{ Q_j B_{s+j-m}^m f_j}^2} 
    \leq
    C 2^{mn} 2^{-s\alpha/2} \ssqnorm[2]{\sum_j\abs{f_j}^2},
\end{equation}
which is the desired unweighted $L^2$ estimate.

To find a crude weighted estimate, write the kernel of $B_{s+j-m}^m$ in
the form
\[
    \int_{2^{s+j-m-1}}^{2^{s+j-m}} \frac1{t^n} \kar_{tS_m}(y)
        \frac{dt}t h(\abs{y}) \sgn\Omega(y).
\]
Since $\abs{y}\leq2^m$ on $S_m$, the kernel is supported in a ball of radius
$2^{s+j}$ centered at the origin, and is bounded by
$2^{n(m-s-j)} \abs{h(\abs{y})}$. Thus
$
    \abs{B_{s+j-m}^mf(x)} \leq c_n 2^{mn} M_{\abs h}f(x)
$,
where, with slight abuse of notation, $M_{\abs h}$ is the maximal
operator of \eqref{e:maxdef} with the radial function
$H(x,y)=\abs{h(\abs{y})}$ in the kernel.  Lemma \ref{l:hstarint} shows
$\abs{h(\abs{\cdot})}$ satisfies condition~\eqref{e:hrect}.

The operator $Q_j$ is convolution by $2^{-nj} \psi(2^{-j}\cdot)$, so
$\abs{Q_j g(x)} \leq C Mg(x)$. Thus
$
    \abs{Q_j B_{s+j-m}^mf(x)} \leq C 2^{mn} M M_{\abs h} f(x)
$.
Let $1<p_1\leq2$ and $v_1\in A_{p_1}$ be arbitrary.
Then
\[
    \bsqnorm[{p_1,v_1}]{\sum_j \abs{ Q_j B_{s+j-m}^m f_j}^2}
    \leq
    C 2^{mn} \bsqnorm[{p_1,v_1}]{\sum_j \abs{ M M_{\abs h} f_j}^2},
\]
and by applying Theorems \ref{t:AndJohn} and \ref{t:maxvect} with
$1<p_1\leq q=2$
we get
\begin{equation}\mylabel{e:badwpest}
    \bsqnorm[{p_1,v_1}]{\sum_j \abs{ Q_j B_{s+j-m}^m f_j}^2} 
    \leq
    C_{p,w} 2^{mn} \ssqnorm[{p_1,v_1}]{\sum_j\abs{f_j}^2}, 
\end{equation}
when $h\in\setH(\sigma)$ and $\sigma$ is large enough.

We treat $Q_j B_{s+j-m}^m$, for fixed $s$ and $m$, as a vector-valued
operator taking values in $\ell^2(\Z)$ and use interpolation with
change of measure between \eqref{e:good2est} and
\eqref{e:badwpest}:

\begin{theorem}[Stein and Weiss \cite{Stei58}]\mylabel{t:SWinter}
For $p_0,$ $p_1\geq1$ suppose that $T$ is a linear operator \textup(possibly
vector-valued\textup) satisfying
$\snorm{Tf}{p_i,v_i} \leq K_i \snorm{f}{p_i,v_i}$ for all
$f\in L^{p_i}_{v_i}$, $i=0,1$. For $0\leq t\leq1$ let
$1/p_t = (1-t)/p_0 + t/p_1$ and $\rho_t=tp_t/p_1$. Define
$v_t = v_0^{1-\rho_t} v_1^{\rho_t}$. Then
$\snorm{Tf}{p_t,v_t} \leq K_0^{1-t} K_1^t \snorm{f}{p_t,v_t}$ for all
$f\in L^{p_t}_{v_t}$.
\end{theorem}

Choose $p_0=2$ and $v_0\equiv1$. Recall that \eqref{e:ThcondC} implies
$w\in A_p$. The $A_p$ properties of $w$ (namely, $w^{1+\alpha}\in
A_{p-\beta}$ for some $\alpha,\beta>0$) allow us to choose $p_1 \in
(1,p]$ such that $v_1=w^{1+\alpha}\in A_{p_1}$ and when $p_t=p$ we have
$\rho_t = 1/(1+\alpha)$. Then, in particular, $v_t=w$ (see
\cite{Duoa86,Wats97}). This yields
\[
    \bsqnorm{\sum_j \abs{ Q_j B_{s+j-m}^m f_j}^2}
    \leq
    C_{p,w} 2^{mn} 2^{-\eta s} \ssqnorm{\sum_j\abs{f_j}^2}
\]
for some $\eta>0$. Thus by \eqref{e:LP2}, the triangle inequality, and
the above estimate, we obtain
\begin{allowdisplaybreaks}
\begin{align*}
    \snorm{\II}{p,w} &= \biggl\| \sum_j \sum_{m\geq0} \sum_{s>Nm}
        Q_j B_{s+j-m}^m Q_jf \biggr\|_{p,w} \\
    &\leq
    C_{p,w} \bsqnorm{ \sum_j \biggl| \sum_{m\geq0} \sum_{s>Nm}
        Q_j B_{s+j-m}^m Q_jf \biggr|^2} \\
    &\leq
    C_{p,w} \sum_{m\geq0} \sum_{s>Nm} \bsqnorm{ \sum_j 
        \abs{Q_j B_{s+j-m}^m Q_jf }^2} \\
    &\leq
    C_{p,w} \sum_{m\geq0} 2^{mn} \sum_{s>Nm} 2^{-\eta s}
        \bsqnorm{ \sum_j \abs{ Q_jf }^2}.
\end{align*}
Note that $\sum_{m\geq0} 2^{mn} \sum_{s>Nm} 2^{-\eta s} \leq c_\eta
\sum_{m\geq0} 2^{mn-\eta Nm} \leq c_{\eta,N} < \infty$, provided that
$N$ is chosen such that $\eta N>n$. We fix such an $N$. Then from above
and by~\eqref{e:LP1},
$
    \snorm{\II}{p,w} \leq
        C_{p,w} \ssqnorm{\sum_j \abs{Q_j f}^2}
        \leq C_{p,w} \snorm{f}{p,w}
$,
which shows $\II$ defines a bounded operator.
\end{allowdisplaybreaks}

\subsection{Term $\III$}

We write $\III$ in the form $\sum_{j\in\Z} \sum_{s<0} \sum_{m\geq0} Q_j B_{s+j-m}^m Q_j f$.
To prove the boundedness of this expression the strategy is to first study the
kernel of $\sum_{s<0} \sum_{m\geq0} Q_j B_{s+j-m}^m$ and show that it is
bounded in absolute value by $2^{-nj} \phi(2^{-j}x)$ for some positive
test function~$\phi$. This gives the estimate
\[
    \biggl|\sum_{s<0} \sum_{m\geq0} Q_j B_{s+j-m}^m Q_j f\biggr| 
    \leq C M(Q_j f)
\]
for some $C>0$ independent of $j$.
Then by \eqref{e:LP2}, Theorem \ref{t:AndJohn}, and finally
estimate~\eqref{e:LP1},
\begin{allowdisplaybreaks}
\begin{align*}
    \snorm{\III}{p,w} &\leq
    C_{p,w} \bsqnorm{\sum_j \biggl| \sum_{s<0} \sum_{m\geq0}
        Q_j B_{s+j-m}^m Q_j f \biggr|^2} \\
    &\leq C_{p,w} \bsqnorm{\sum_j \abs{ M(Q_j f)}^2} \\
    &\leq C_{p,w} \bsqnorm{\sum_j \abs{ Q_j f}^2}
    \leq C_{p,w} \snorm{f}{p,w}. 
\end{align*}
\end{allowdisplaybreaks}

The kernel of the operator $Q_j B_{s+j-m}^m$ is the convolution
$\psi_j * b_{s+j-m}^m(x)$, so it is equal to
\begin{multline*}
    \int_\Rn 2^{-nj} \psi(2^{-j}(x-y))
    \int_{2^{s+j-m-1}}^{2^{s+j-m}} \frac1{t^n} \kar_{tS_m}(y) \frac{dt}t\;
    h(\abs{y}) \sgn\Omega(y)\,dy \\
    =
    K_{s,j,m}(x) + 2^{-nj} \psi(2^{-j}x) L_{s,j,m},
\end{multline*}
where
\begin{align*}
  K_{s,j,m}(x) &=
    \int_\Rn \biggl[
    2^{-nj} (\psi(2^{-j}(x-y)) - \psi(2^{-j}x))   \\
    &\hspace*{5em}\ \times
    \int_{2^{s+j-m-1}}^{2^{s+j-m}} \frac1{t^n} \kar_{tS_m}(y) \frac{dt}t\;
    h(\abs{y}) \sgn\Omega(y)
    \biggr]\,dy, \displaybreak[0]\\
  L_{s,j,m} &=
    \int_\Rn 
    \int_{2^{s+j-m-1}}^{2^{s+j-m}} \frac1{t^n} \kar_{tS_m}(y) \frac{dt}t\;
    h(\abs{y}) \sgn\Omega(y)\,dy. 
\end{align*}

To estimate $\sum_{s<0} \sum_{m\geq0} K_{s,j,m}(x)$, define
\[
    k_{s+j}(y) = \sum_{m\geq0}
    \int_{2^{s+j-m-1}}^{2^{s+j-m}} \frac1{t^n} \kar_{tS_m}(y) \frac{dt}t\;
    h(\abs{y}) \sgn\Omega(y).
\]
Then
\begin{equation}\mylabel{e:KsjmA} 
    \biggl| \sum_{s<0} \sum_{m\geq0} K_{s,j,m}(x) \biggr|
    \leq
    \sum_{s<0} \int_\Rn 2^{-nj} \abs{\psi(2^{-j}(x-y)) - \psi(2^{-j}x)}
        \abs{k_{s+j}(y)} \, dy. 
\end{equation}
Note that $k_{s+j}(y)\not=0$ if and only if $y\in tS_m$ for some
$m\geq0$ and
$t$ in the interval $(2^{s+j-m-1}, 2^{s+j-m})$. Since
$S_m\subset B(0,2^m)$, the support of 
$k_{s+j}$ is contained in the ball $B(0,2^{s+j})$. Hence in the integral in
\eqref{e:KsjmA} we have $2^{-j}\abs{y}\leq 2^s < 1$. 
Therefore there exists a positive $\phi_0 \in \Sch$ such that
\[
    2^{-nj} \abs{\psi(2^{-j}(x-y)) - \psi(2^{-j}x)}
    \leq
    2^s 2^{-nj} \phi_0(2^{-j}x)
\]
(use the mean value theorem and
let $\phi_0(x)$ majorize $\sup_{\abs{y}<1} \abs{\nabla \psi(x-y)}$).

This allows us to estimate \eqref{e:KsjmA} further. We get
\[
    \biggl| \sum_{s<0} \sum_{m\geq0} K_{s,j,m}(x) \biggr|
    \leq
    2^{-nj} \phi_0(2^{-j}x) \sum_{s<0} 2^s \int_{\Rn}
        \abs{k_{s+j}(y)} \, dy, 
\]
and
\[
    \int_\Rn \abs{k_{s+j}(y)} \, dy
    \leq
    \sum_{m\geq0} \int_{2^{s+j-m-1}}^{2^{s+j-m}} \frac1{t^n}
    \int_{tS_m} \abs{h(\abs{y})} \, dy \frac{dt}t.
\]
According to Lemma \ref{l:hstarint} the innermost integral is bounded
by $Ct^n\abs{S_m}$. Hence
$\int_\Rn \abs{k_{s+j}(y)} \, dy \leq
    C \sum_{m\geq0} \abs{S_m} = C \snorm{\Omega}{1}$.
Thus we get that
\begin{equation}\mylabel{e:Ksjmest}
    \biggl| \sum_{s<0} \sum_{m\geq0} K_{s,j,m}(x) \biggr|
    \leq
    C \snorm{\Omega}{1} 2^{-nj} \phi_0(2^{-j}x).
\end{equation}

To handle $\sum_{s<0} \sum_{m\geq0} L_{s,j,m}$ we first
show absolute convergence. We have 
$
    \abs{L_{s,j,m}(x)}
	\leq
    \int_{2^{s+j-m-1}}^{2^{s+j-m}} \int_{S_m}  \abs{h(t\abs{y})} \,dy \,dt/t
$.
Summation over $s<0$ results in
$
    \int_0^{2^{j-m-1}} \int_{S_m}  \abs{h(t\abs{y})} \,dy \,  \frac{dt}t,
$
which is bounded by
$
    \int_0^{2^{j}} \int_{S_m}  \abs{h(t\abs{y})} \,dy \,  \frac{dt}t.
$
The sum over $m\geq0$ results in the bound
$
    \int_0^{2^{j}} \int_{S}  \abs{h(t\abs{y})} \,dy \,  \frac{dt}t.
$
Lemma \ref{l:hSint} shows the $t$-integral over $(0,1)$
is bounded by
$C(h,n) \snorm{\Omega}{\LlogL}$.
For $j>0$ the remaining part, over $(1,2^j)$, is clearly bounded by
$C(h,n)\max\{j,1\} \* \snorm{\Omega}{1}$ by Lemma \ref{l:hstarint}.
What we need, however, is a bound independent of~$j$ and the support of $h$.

The above argument allows us to change the order of summation and
integration. By also changing variables in both integrals we get
\begin{multline}\mylabel{e:LsjmA}
    \sum_{s<0} \sum_{m\geq0} L_{s,j,m} =
    \sum_{m\geq0} \int_{S_m} \int_0^{2^{j-m-1}\abs{y}} h(t) \frac{dt}t
        \sgn\Omega(y)\,dy \\
    =
    \sum_{m\geq0} \int_{S_m} \int_{2^j\abs{y}}^{2^{j-m-1}\abs{y}} h(t) \frac{dt}t
        \sgn\Omega(y)\,dy
    +
    \int_{S} \int_0^{2^j\abs{y}} h(t) \frac{dt}t
        \sgn\Omega(y)\,dy.
\end{multline}
Using $\int_S \sgn\Omega(y)\,dy=0$ and the fact that
$\int_0^{2^j} h(t) \,dt/t$ is independent of $y$,
the last term is equal to
$
    \int_S \int_{2^j}^{2^j\abs{y}} h(t) \frac{dt}t
        \sgn\Omega(y)\,dy,
$
which is bounded by
\begin{equation}\mylabel{e:LsjmB}
    \int_S \biggl| \int_{2^j}^{2^j\abs{y}} \abs{h(t)} \frac{dt}t \biggr|\,dy
    \leq
    C \int_S \big| \log\abs{y} \big|\,dy
    \leq
    C (\snorm{\Omega}{\LlogL} + 1), 
\end{equation}
the two inequalities following from Lemma \ref{l:Hlog} and 
from \eqref{e:Slogint}.

Similarly by Lemma \ref{l:Hlog} and \eqref{e:mSmsum} of Lemma \ref{l:Sints},
\begin{align}\mylabel{e:LsjmC}
     \sum_{m\geq0} \int_{S_m} \int_{2^{j-m-1}\abs{y}}^{2^j\abs{y}} \abs{h(t)} \frac{dt}t
        \,dy
     &\leq
     C \sum_{m\geq0} \abs{S_m} \log 2^{m+1} \\
     &\leq
     C \snorm{\Omega}{\LlogL}. \notag
\end{align}

Finally, putting together \eqref{e:LsjmA}, \eqref{e:LsjmB}, and
\eqref{e:LsjmC} we get
\begin{equation}\mylabel{e:Lsjmest}
    \biggl|\sum_{s<0} \sum_{m\geq0} L_{s,j,m} \biggr|
    \leq
    C  (\snorm{\Omega}{\LlogL}+1).
\end{equation}
Note that this estimate is independent of $j$ and the support of $h$.

Combining \eqref{e:Ksjmest} and \eqref{e:Lsjmest} we see that the
absolute value of the kernel of
$\sum_{s<0} \sum_{m\geq0} Q_j B_{s+j-m}^m$ is bounded by
$
    C 2^{-nj} (\phi_0(2^{-j}x) + \abs{\psi(2^{-j}x)}),
$
which in turn is bounded by $2^{-nj} \phi(2^{-j}x)$ for some positive
$\phi\in\Sch$. 
Now the argument given at the beginning of this section shows that
$\III$ is bounded on $\Lpw$.

This completes the proof of Theorem \ref{t:Thbded}.

\begin{proof}[Proof of theorem \ref{t:Thinfbded}]
We show that when $p=2$ and $h\in\Li(\Rp)$ we can take $r=1$ in the
above proof.

When $p=2$ to show the boundedness of term $\I$ it is enough to show
only \eqref{e:Gjpw}. Estimate \eqref{e:Gjsuppw} and the interpolation
argument between \eqref{e:Gjpw} and \eqref{e:Gjsuppw} is not needed.
When $h\in\Li(\Rp)$ the application of H\"older's inequality in
\eqref{e:anotherholder} works with $r=1$ and $r'=\infty$, with the usual
meaning for the first factor on the right-hand side of
\eqref{e:anotherholder}. Hence condition~\ref{e:ThcondB} is enough.

The only place where the stronger condition (with $r>1$) on the weight
is used in the proof for $\II$ is when Theorem~\ref{t:MSh} is applied to
the maximal operator $M_{\abs{h}}$. But when $h\in\Li(\Rp)$ we have
$M_{\abs{h}}f \leq \snorm{h}{\infty} Mf$, where $M$ is the
Hardy-Littlewood maximal function, hence the estimates for $\II$ now
require only $w\in A_2$.

Finally, the estimates for term $\III$ use only the fact that $w$ is
in $A_2$.
\end{proof}
 

\appendix
\section{Other representation formulas}\mylabel{s:hrep}

The representation formula given in Theorem~\ref{t:srep} is only for
convolution operators. This is not essential. For non-convolution
type operators we have the following result:

\begin{theorem}\mylabel{t:srepk}
Assume $\Omega\in L^1(\Sn)$ is positively homogeneous of degree zero, ${k\in
\Li(\Rn\times\Rn)}$,
and either $f\in\Li(\Rn)$ with compact support
or $f\in L^p(\Rn)$, ${1\leq p<\infty}$. For $\epsilon>0$ define the operator
\[
    T^{(k)}_\epsilon f(x) = \int_{|x-y|>\epsilon} \frac{\Omega(x-y)}{|x-y|^n}
    k(x,y) f(y)\,dy
\]
and let
\begin{subequations}\mylabel{e:srepk}
\begin{equation}\mylabel{e:srepkX}
    A_{\epsilon,t} f(x) = \frac1{t^n} \int_{tS\setminus 
    B(0,\epsilon)} f(x-y)k(x,x-y)\sgn\Omega(y)\,dy,
\end{equation}
where $S=S_\Omega$ is the star-shaped set associated with $\Omega$.
Then for almost all $x\in\Rn$ the representation formula
\begin{equation}\mylabel{e:srepkY}
    T^{(k)}_\epsilon f(x) = n \int_0^\infty A_{\epsilon,t} f(x) \frac{dt}t
\end{equation}
\end{subequations}
holds and the integrals in \eqref{e:srepkX} and \eqref{e:srepkY} converge 
absolutely.
\end{theorem}

The proof is practically the same as the proof of
Theorem~\ref{t:srep} given in Section~\ref{s:srep} and is therefore omitted.
Proving Theorem~\ref{t:srepk} is actually easier, since $k(x,y)$ is a
bounded function as compared to $h(\abs{x-y})$, which maybe unbounded.

For principal value operators we can derive similar formulas. 
E.g., for convolution operators we have the following result
(recall from remark~\ref{r:limTheps} that the principal value operator
is well defined on test functions).

\begin{theorem}\mylabel{t:hrep}
Suppose $\Omega\in\LlogLS$ is positively homogeneous of degree~$0$,
$\int_\Sn\Omega(\theta)\,d\theta=0$, $h\in\setH(1)$, and
\begin{equation}\mylabel{e:hdini}
    \int_0^1 \abs{h(t)-h(0)}\frac{dt}t < \infty.
\end{equation}
Let
\begin{equation}\mylabel{e:Thpv}
    Tf(x) = \pv \int_\Rn \frac{\Omega(x-y)}{\abs{x-y}} h(\abs{x-y}) f(y) \,dy
\end{equation}
and
\begin{subequations}\label{e:Thpvrepboth}
\begin{equation}\mylabel{e:Athpv}
    A_tf(x) = \frac1{t^n} \int_{tS} h(\abs{y})f(x-y)\sgn\Omega(y)\,dy,
\end{equation}
where $S=S_\Omega$ is the star-shaped set associated with $\Omega$ and 
$f\in C_0^1(\Rn)$. Then the integral
$\int_0^\infty A_tf(x) t^{-1}\,dt$ converges absolutely and
\begin{equation}\mylabel{e:Thpvrep}
    Tf(x) = n \int_0^\infty A_t f(x) \frac{dt}t + h(0)c_\Omega f(x)
\end{equation}
\end{subequations}
for all $x\in\Rn$ and all $f\in C_0^1(\Rn)$, where $c_\Omega =
\frac1n \int_\Sn \Omega(\theta)
\log\abs{\Omega(\theta)}\,d\theta$. 
\end{theorem}

Our proof of the above theorem requires the Dini-condition
\eqref{e:hdini}, even though such a condition is not necessary for any
of the boundedness or convergence results discussed in
Section~\ref{s:res}.  It remains open whether there is a result
similar to Theorem~\ref{t:hrep} but without condition
\eqref{e:hdini}. 

Note the extra term in \eqref{e:Thpvrep} when compared to the
representation formula \eqref{e:srepY} for truncated operators. Using
the pointwise convergence of the truncated operators (remark
\ref{r:limTheps}) we have that
$
   \int_0^\infty A_{\epsilon,t} f(x)\,dt/t
$
converges to
$
    \int_0^\infty A_t f(x)\,dt/t + \frac1n h(0) c_\Omega f(x)
$
as $\epsilon\tozerop$.

As above there are corresponding results for non-convolution type
operators. We discuss the Calder\'on commutators as an example:
Given a function $a(x)$ with $\nabla a\in\Li(\Rn)$
the $k$th Calder\'on commutator is
$
    C^{(a,k)} f(x) = \lim_{\epsilon\tozerop} 
    C^{(a,k)}_\epsilon f(x)
$,
$k = 1$, $2$, \ldots,
where the truncated operators are defined by
\[
    C^{(a,k)}_\epsilon f(x) = \int_{\abs{x-y}>\epsilon}
    \frac{\Omega(x-y)}{\abs{x-y}^n} \left( \frac{a(x)-a(y)}{\abs{x-y}} \right)^k
    f(y) \, dy.
\]
For the truncated operators we have the representation formula
\[
    C^{(a,k)}_\epsilon f(x) = 
      n \int_0^\infty \frac1{t^n} \int_{tS\setminus B(0,\epsilon)}
      f(x-y)  \left(\frac{a(x)-a(x-y)}{\abs{y}}\right)^k
      \sgn\Omega(y)\,dy\frac{dt}t,
\]
when $\Omega$ and $f$ are as in Theorem~\ref{t:srepk}. This follows from
Theorem~\ref{t:srepk}.

If $a$ satisfies an additional regularity condition we get a
representation formula for the principal value operator. Let
\[
    w_x(t) = \sup_{\theta\in\Sn} 
    \frac{\abs{a(x)-a(x-t\theta)-\nabla a(x)\cdot t\theta}}{t}
\]
and assume that $\int_0^1 w_x(t)t^{-1} \,dt$ is finite for each
$x\in\Rn$. This corresponds to condition \eqref{e:hdini} for $h$ in
Theorem~\ref{t:hrep}. 
If $\Omega \in \LlogLS$ is homogeneous of degree zero
and satisfies
$
    \int_\Sn\theta^\alpha\Omega(\theta)\,d\theta=0
$
for all multi-indices $\alpha$
with $\abs{\alpha}=k$, then
\begin{multline*}
    C^{(a,k)}f(x) =
    c_\Omega \int_\Sn \int_1^{\rho(\theta)} 
        (\nabla a(x) \cdot \theta)^k
        \Omega(\theta)\frac{dr}r\,d\theta \, f(x) \\
    + n \int_0^\infty \frac1{t^n} \int_{tS}
      f(x-y)  \left(\frac{a(x)-a(x-y)}{\abs{y}}\right)^k
      \sgn\Omega(y)\,dy\frac{dt}t
\end{multline*}
with $f$ and $c_\Omega$ as in Theorem \ref{t:hrep}. The proof is almost
identical with the proof of Theorem~\ref{t:hrep} given below and is
therefore omitted. Note that when $\Omega \in \LlogLS$ the multiplier
$
    c_\Omega \int_\Sn \int_1^{\rho(\theta)} (\nabla a(x) \cdot 
    \theta)^k\Omega(\theta)\,dr/r\,d\theta
$
is a bounded function of $x$.

\subsection{Proof of Theorem \ref{t:hrep}}

The proof is a generalization of an argument for the case $h\equiv1$
that was given in a preprint version of \cite{Wats97}.

\subsubsection{Absolute convergence}

We begin the proof by showing $\int_0^\infty A_tf(x) t^{-1}\,dt$
converges absolutely. 
%
%
%
%
When $0\leq t\leq1$ we use
$\int_S \sgn\Omega(y)\,dy = \frac1n\int_\Sn\Omega(\theta)\,d\theta=0$
to bound $\int_0^\infty \abs{A_tf(x)} t^{-1}\,dt$  by
\begin{multline}\mylabel{e:hrepA}
   \int_0^1 \int_S \abs{h(t\abs{y}) (f(x-ty)-f(x))} \,dy
                \frac{dt}t \\
    + \abs{f(x)} \int_0^1 \int_S \abs{h(t\abs{y})-h(0)} \,dy 
                \frac{dt}t 
    + \int_1^\infty \int_S \abs{h(t\abs{y}) f(x-ty)}\,dy
                \frac{dt}t.
\end{multline}

To show the first term is finite, we use the fact that since
$\abs{f(x-ty)-f(x)}$ is bounded by both $2\snorm{f}{\infty}$ and
$\snorm{\nabla f}{\infty} t\abs{y}$, it is bounded by $C_f
t\abs{y}/(1+t\abs{y})$. Making a change of variables in the $t$-integral
($t\abs{y}$ replaced with $t$) we get that the first term of
\eqref{e:hrepA} is at most a constant times
%
%
%
%
\[
    \int_S \int_0^{\abs{y}} \abs{h(t)} \frac{1}{1+t} \,dt \,dy
    \leq
    \int_S \int_0^1 \abs{h(t)} \,dt\,dy + 
        \int_S \int_0^{\bmax{\abs{y}}{1}} \abs{h(t)} \,\frac{dt}t \,dy,
\]
where $\bmax{\abs{y}}{1}=\max\{\abs{y},1\}$. Using Lemma \ref{l:Hlog} in
the last term shows that this is bounded by
$
    C\int_S (1+\logp\abs{y}) \,dy
    \leq
    C\snorm{\Omega}{\LlogL}.
$

In the second term of \eqref{e:hrepA}
we change the order of integration and then change variables
in the $t$-integral as above. This gives,
without the factor $\abs{f(x)}$, the estimate
\[
    \int_S \int_0^1 \abs{h(t)-h(0)} \frac{dt}t \,dy
        + \int_S \int_1^\bmax{\abs{y}}{1} \abs{h(t)-h(0)} \frac{dt}t \,dy.
\]
By \eqref{e:hdini} and \eqref{e:Sint} of Lemma
\ref{l:Sints} the first term of this expression is at most $C
\snorm{\Omega}{1}$. The second term is bounded by
\[
    \int_S \int_1^\bmax{\abs{y}}{1} \abs{h(t)}\,\frac{dt}t\,dy +
        \abs{h(0)} \int_S \logp\abs{y}\,dy,
\]
and by Lemma \ref{l:Hlog} this is at most
$(C+\abs{h(0)})\snorm{\Omega}{\LlogL}<\infty$. 

Let $R>0$ be such that $\supp f\subset B(0,R)$. If $x-z\in\supp f$ then
$\abs{z}\leq\abs{r\theta-x}+x\leq R+\abs{x}$.
The change of coordinates $z=ty$ shows the third term of \eqref{e:hrepA} is
\[
    \int_1^\infty \int_{tS} \abs{h(\abs{z}) f(x-z)} \,dz\frac{dt}{t^{n+1}}
    \leq
    \snorm{f}{\infty}
    \int_1^\infty \frac{dt}{t^{n+1}}
    \int_{B(0,R+\abs{x})} \abs{h(\abs{z})} \,dz.
\]
Lemma \ref{l:hstarint} implies this is bounded by $C\snorm{f}{\infty}
(R+\abs{x})^n$.

\subsubsection{Representation formula}

To prove the representation formula write
\begin{equation}\mylabel{e:hrepB}
    \int_0^\infty A_t f(x)\frac{dt}t = \I + \II + f(x) \III,
\end{equation}
where, corresponding to \eqref{e:hrepA},
\begin{allowdisplaybreaks}
\begin{align*}
    \I &= \int_0^1 \int_S h(t\abs{y}) (f(x-ty)-f(x)) \sgn\Omega(y)\,dy\frac{dt}t,\\
   \II &= \int_1^\infty \int_S h(t\abs{y}) f(x-ty) \sgn\Omega(y)\,dy\frac{dt}t,\\
  \III &= \int_0^1 \int_S (h(t\abs{y})-h(0)) \sgn\Omega(y)\,dy\frac{dt}t.
\end{align*}
The above computations also show that the multiple integrals defining
each of the terms $\I$, $\II$ and $\III$ converge absolutely.
\end{allowdisplaybreaks}

Making the change of variables $\eta=ty$ in the $y$-integral of $\I$ 
%
%
%
%
and using polar coordinates $\eta=r\theta$ we get
\[
    \I = \int_0^1 \int_\Sn \int_0^{t\rho(\theta)} h(r) 
        (f(x-r\theta) - f(x)) \sgn\Omega(\theta) r^{n-1}\,dr
        \,d\theta \frac{dt}{t^{n+1}}.
\]
Changing the order of integration to make the $t$-integral the inner
most (see figure~\ref{f:Idom}) gives
\begin{figure}
\begin{center}
{\setlength{\unitlength}{3.5cm}
\begin{picture}(1.3,1.15)(-0.15,-0.1)
    \put(-0.05,0){\vector(1,0){1.2}}
    \put(0,-0.05){\vector(0,1){1.1}}
    \put(0,0){\line(4,3){1.1}}
    {\thicklines\put(0,0){\line(4,3){1}}
    \put(0,0){\line(1,0){1}}
    \put(1,0){\line(0,1){0.75}}}
    \put(1,-0.025){\line(0,1){0.05}}	
    \put(-0.025,0.75){\line(1,0){0.05}}	
    \multiput(0.05,0.75)(0.0545,0){18}{\line(1,0){0.025}}
    \put(1,-0.05){\makebox(0,0)[t]{$1$}}
    \put(1.125,-0.05){\makebox(0,0)[t]{$t$}}
    \put(-0.025,0.75){\makebox(0,0)[r]{$\rho(\theta)$}}
    \put(-0.025,1){\makebox(0,0)[r]{$r$}}
    \put(0.25,0.3){\rotatebox{36.9}{$r=t\rho(\theta)$}}
    \put(0.6,0.15){\makebox(0,0)[b]{$\I$, $\III$}}
\end{picture}}
\qquad
\qquad
{\setlength{\unitlength}{3.5cm}
\begin{picture}(1.3,1.15)(-0.15,-0.1)
    \put(-0.05,0){\vector(1,0){1.2}}
    \put(0,-0.05){\vector(0,1){1.1}}
    \put(0,0){\line(4,3){1.1}}
    {\thicklines
        \put(0.3,0){\line(0,1){0.225}}
        \put(0.3,0){\line(1,0){0.8}}
        \put(0.3,0.225){\line(4,3){0.8}}}
    \put(0.3,-0.025){\line(0,1){0.05}}	
    \put(-0.025,0.225){\line(1,0){0.05}}	
    \multiput(0.3,0.225)(0.054,0){15}{\line(1,0){0.025}}
    \put(0.3,-0.05){\makebox(0,0)[t]{$1$}}
    \put(1.125,-0.05){\makebox(0,0)[t]{$t$}}
    \put(-0.025,0.225){\makebox(0,0)[r]{$\rho(\theta)$}}
    \put(-0.025,1){\makebox(0,0)[r]{$r$}}
    \put(0.25,0.3){\rotatebox{36.9}{$r=t\rho(\theta)$}}
    \put(0.6,0.08){\makebox(0,0)[b]{$\II_2$}}
    \put(0.8,0.35){\makebox(0,0)[b]{$\II_1$}}
\end{picture}}
\end{center}
\caption{Domain of integration in $\I$ and $\III$ (left) and
         in $\II$ (right).}\mylabel{f:IIdom}
\mylabel{f:Idom}
\end{figure}
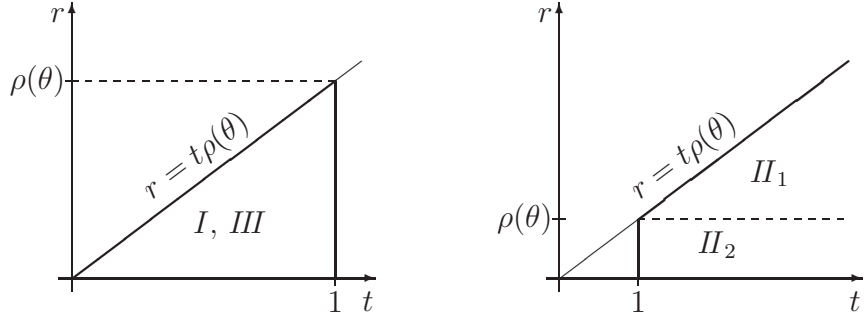
\[
    \I = \int_\Sn \int_0^{\rho(\theta)} h(r) 
        (f(x-r\theta) - f(x)) \sgn\Omega(\theta)
        \int_{r/\rho(\theta)}^1 \frac{dt}{t^{n+1}}
        r^{n-1}\,dr \,d\theta.
\]
Now use the identity
$
    \int_{r/\rho(\theta)}^1 t^{-n-1} \,dt
    = n^{-1} \bigl( (\rho(\theta)/r)^n - 1 \bigr)
$
%
%
%
%
to get $\I = \I_1 + \I_2$, where
\begin{allowdisplaybreaks}
\begin{align*}
    \I_1 &= \frac1n \int_\Sn \int_0^{\rho(\theta)} h(r) (f(x-r\theta) - f(x))
        \Omega(\theta) \frac{dr}r \,d\theta, \\
    \I_2 &= - \frac1n \int_\Sn \int_0^{\rho(\theta)} h(r) (f(x-r\theta) - f(x))
        \sgn\Omega(\theta) r^{n-1}\,dr \,d\theta.
\end{align*}
The%
\end{allowdisplaybreaks}
multiple integral in $\I_2$ converges absolutely since the
corresponding integral with absolute values in the integrand can be
bounded by
$    C_h \snorm{f}{\infty} \snorm{\Omega}{1}
$
%
%
by Lemma \ref{l:hstarint}. We have already
shown the integral in $\I$ is absolutely convergent, so by the triangle
inequality also the multiple integral in $\I_1$ converges absolutely.

We do similar computations for term $\II$:
\[
    \II = \int_1^\infty \int_\Sn \int_0^{t\rho(\theta)} h(r) f(x-r\theta)
        \sgn\Omega(\theta) r^{n-1}\,dr\,d\theta\frac{dt}{t^{n+1}}
\]
%
%
%
%
and changing the order of integration, see figure~\ref{f:IIdom}, gives
$\II = \II_1 + \II_2$, where
\begin{allowdisplaybreaks}
\begin{align*}
    \II_1 &= \frac1n \int_\Sn \int_{\rho(\theta)}^\infty h(r) f(x-r\theta)
            \Omega(\theta) \frac{dr}r\,d\theta, \\
    \II_2 &= \frac1n \int_\Sn \int_0^{\rho(\theta)} h(r) f(x-r\theta)
            \sgn\Omega(\theta) \frac{dr}r\,d\theta.
\end{align*}
The first term corresponds to the triangular region and the
second term to the rectangle.
Since they are obtained by decomposing the
domain of integration of the absolutely convergent integral in $\II$ 
into two disjoint sets, the multiple integrals in $\II_1$ and $\II_2$ 
converge absolutely.

\end{allowdisplaybreaks}

Finally we get $\III = \III_1 + \III_2$, where
\begin{allowdisplaybreaks}
\begin{align*}
    \III_1 &= \frac1n \int_\Sn \int_0^{\rho(\theta)} (h(r)-h(0))
            \Omega(\theta) \frac{dr}r\,d\theta, \\
    \III_2 &= -\frac1n \int_\Sn \int_0^{\rho(\theta)} (h(r)-h(0))
            \sgn\Omega(\theta) r^{n-1}\,dr\,d\theta
\end{align*}
(this is similar to~$\I$). 
Absolute convergence of all integrals is again easy to show.

\end{allowdisplaybreaks}

Now notice that
\begin{equation}\mylabel{e:hrepH}
    \I_2 + \II_2 + f(x)\III_2
    =
    \frac1n f(x) h(0) \int_\Sn \Omega(\theta) \,d\theta = 0
\end{equation}
and write
\begin{multline}\mylabel{e:hrepG}
    \I_1 + \II_1
    = \frac1n \int_\Sn \left( \int_0^1 + \int_1^{\rho(\theta)} \right)
            h(r) (f(x-r\theta)-f(x)) \Omega(\theta) \frac{dr}r\,d\theta \\
    + \frac1n \int_\Sn \left( \int_1^\infty + \int_{\rho(\theta)}^1\right)
            h(r) f(x-r\theta) \Omega(\theta) \frac{dr}r\,d\theta
    = A + B + C + D.
\end{multline}
It is easy to see all four of the above multiple integrals are finite (actually
they converge absolutely).

We get from \eqref{e:hrepG} and the definition of $\III_1$ that
\begin{multline}\mylabel{e:hrepC}
    \I_1 + \II_1 + f(x)\III_1 
    =
    (A+C) + \bigl[ (B+D) + f(x)\III_1 \bigr] \\        
    =
        \frac1n \int_\Sn \int_0^\infty h(r)
            \bigl(f(x-r\theta)-f(x)\kar_{(0,1)}(r)\bigr) 
            \Omega(\theta) \frac{dr}r\,d\theta \\
       + \frac1n f(x) \Biggl[ - \int_\Sn \int_1^{\rho(\theta)}
            h(r) \Omega(\theta) \frac{dr}r\,d\theta  \\
            + \int_\Sn \int_0^{\rho(\theta)} (h(r)-h(0))
                \Omega(\theta) \frac{dr}r\,d\theta \Biggr].
\end{multline}
The expression inside the square brackets is equal to
\begin{equation}\mylabel{e:hrepD}
    \int_\Sn \int_0^1 (h(r)-h(0)) \Omega(\theta)\frac{dr}r \,d\theta
    -
    h(0) \int_\Sn \int_1^{\rho(\theta)} \Omega(\theta)\frac{dr}r \,d\theta
    =
    -h(0)c_\Omega,
\end{equation}
since $\int_\Sn\Omega(\theta)\,d\theta=0$ implies the first term is zero.

On the other hand when $f\in C_0^1(\Rn)$,
\[
    Tf(x) = \int_\Sn \int_0^\infty
                h(r) (f(x-r\theta)-f(x)\kar_{(0,1)}(r))
                    \Omega(\theta) \frac{dr}r \,d\theta,
\]
%
%
%
therefore by \eqref{e:hrepB}, \eqref{e:hrepH}, \eqref{e:hrepC},
and \eqref{e:hrepD}
\begin{align*}
    Tf(x) &
    = n\bigl(\I_1 + \II_1 + f(x)\III_1\bigr) + h(0)c_\Omega f(x)
          \\
          &
    = n \int_0^\infty A_t f(x) \frac{dt}t + h(0)c_\Omega f(x).
\end{align*}
This completes the proof of Theorem \ref{t:hrep}.

 {\footnotesize

\bibliographystyle{amsplain}

\addcontentsline{toc}{section}{References}
\bibliography{mrabbrev,analysis}
}

\medbreak

{\footnotesize
{\bf Address}: Department of Mathematics - Hill Center;
    Rutgers, the State University of New Jersey;
    110~Frelinghuysen Rd; 	
    Piscataway NJ 08854-8019; USA

{\bf E-mail address}: ojanen@math.rutgers.edu
}
\end{document}